\documentclass{amsart}
\usepackage[cp850]{inputenc}
\usepackage{diagrams}%

\RequirePackage{amssymb} 
\RequirePackage{amsmath} 

\usepackage{bookman}%

\newtheorem{theorem}{Theorem}[section]
\newtheorem{proposition}[theorem]{Proposition}
\newtheorem{lemma}[theorem]{Lemma}
\newtheorem{corollary}[theorem]{Corollary}
\newtheorem{problem}[theorem]{Problem}

\theoremstyle{definition}

\newtheorem{example}[theorem]{Example}

\theoremstyle{remark}
\newtheorem{remark}[theorem]{Remark}

\numberwithin{equation}{section}

\DeclareMathOperator{\calM} {\mathcal M}%
\DeclareMathOperator{\calQ} {\mathcal Q}%
\DeclareMathOperator{\calS} {\mathcal S}%
\DeclareMathOperator{\calV} {\mathcal V}%
\DeclareMathOperator{\calF} {\mathcal F}%

\DeclareMathOperator{\content} {\mathbf{c}}%
\DeclareMathOperator{\f} {\boldsymbol{f}}%
\DeclareMathOperator{\F} {\boldsymbol{F}}%
\DeclareMathOperator{\Ima} {Im}%
\DeclareMathOperator{\Kr} {Kr}%
\DeclareMathOperator{\Max} {Max}%
\DeclareMathOperator{\Na} {Na}%
\DeclareMathOperator{\oF} {\boldsymbol{\overline{F}}}%
\DeclareMathOperator{\rad} {rad}%
\DeclareMathOperator{\sdDP} {\dot{\star}^{\it {D_P}} }%
\DeclareMathOperator{\sM} {\star_{\it M}}%
\DeclareMathOperator{\sP} {\star_{\it P}}%
\DeclareMathOperator{\Spec} { Spec}%
\newcommand{\balf}
 {\renewcommand{\theenumi}{(\alph{enumi})}
 \renewcommand{\labelenumi}{\bf\theenumi}
 \begin{enumerate}}
\newcommand{\ealf}
 {\end{enumerate}
 \renewcommand{\theenumi}{\arabic{enumi}}
 \renewcommand{\labelenumi}{\theenumi.}}
\newcommand{\bara}
 {\renewcommand{\theenumi}{(\arabic{enumi})}
 \renewcommand{\labelenumi}{\bf\theenumi}
 \begin{enumerate} }
\newcommand{\eara}
 {\end{enumerate}
 \renewcommand{\theenumi}{\arabic{enumi}}
 \renewcommand{\labelenumi}{\theenumi.}}
\newcommand{\brom}
 {\renewcommand{\theenumi}{(\roman{enumi})}
 \renewcommand{\labelenumi}{\bf\theenumi}
 \begin{enumerate} }
\newcommand{\erom}
 {\end{enumerate}
 \renewcommand{\theenumi}{\arabic{enumi}}
 \renewcommand{\labelenumi}{\theenumi.}}

\begin{document}
\title{Local--Global Properties for Semistar Operations}
%
%
\author{M. Fontana}
\thanks{First author was supported in part by a research grant MIUR 2001/2002 (Cofin 2000-MM01192794)}
\address{Dipartimento di Matematica, Universit{\`a} degli Studi ``Roma Tre''.
\newline
\indent Largo San Leonardo Murialdo, 1. 00146--Roma. ITALY}
\email{fontana@mat.uniroma3.it}%
\urladdr{
http://www.mat.uniroma3.it/dipartimento/membri/fontana\_homepage.html}%
%
%
\author{P. Jara}
\address{Department of Algebra. University of Granada. 18071--Granada. SPAIN}
\email{pjara@ugr.es}%
\urladdr{http://www.ugr.es/local/pjara}%
%
%
\author{E. Santos}
\thanks{The second and third authors were partially supported by %
        DGES BMF2001--2823 and %
        FQM-266 (Junta de Andaluc{\'\i}a Research Group).}
\address{Department of Algebra. University of Granada. 18071--Granada. SPAIN}

\email{esantos@ugr.es}%
%
%
{\subjclass[2000] {Primary 13F05; Secondary 13G05}}
\date{April 28, 2003}
\keywords{Pr{\"u}fer domain, Nagata ring, Kronecker function ring.}

\begin{abstract}
We study the ``local" behavior of several relevant properties concerning semistar
operations, like finite type, stable, spectral, e.a.b. and a.b. We deal with the
``global" problem of building a new semistar operation on a given integral domain, by
``gluing" a given homogeneous family of semistar operations defined on a set of
localizations. We apply these results for studying the local--global behavior of the
semistar Nagata  ring  and the semistar Kronecker function  ring. We prove that an
integral domain $D$ is a Pr{\"u}fer $\star$--multiplication domain if and only if all
its localizations $D_{\!P}$ are Pr{\"u}fer $\star_P$--multiplication domains.
\end{abstract}

\maketitle

\section*{Introduction}

Krull's theory on ideal systems and star operations was motivated by the construction
of Kronecker function rings in a more general context than that of algebraic integers,
originally considered by L. Kronecker. The theory developed by Krull requires some
restrictions on the integral domain, which has to be integrally closed, and on the star
operation, which has to be e.a.b. \cite[Section 32]{GILMER:1972}. Semistar operations,
introduced by Okabe and Matsuda \cite{Okabe/Matsuda:1994}, lead to very general theory
of Kronecker function rings, also in case of non necessarily integrally closed domains
(cf. \cite{Okabe/Matsuda:1997}, \cite{Matsuda:1998}, \cite{Koch:2000}, \cite{FL1},
\cite{FL2} and \cite{FL3}).

Semistar operations are an appropriate tool for extending the theory of Pr{\"u}fer
domains and, more generally, of Pr{\"u}fer $v$--multiplication domains (cf.
\cite{Griffin:1967}, \cite{Mott/Zafrullah:1981} and \cite{Kang:1989}) to the non
necessarily integrally closed domains case. Let $\star$ be a semistar operation of
finite type on an integral domain $D$ (the formal definition is recalled in Section 1),
then $D$ is called a Pr{\"u}fer $\star$--multiplication domain (for short, P$\star$MD)
if each nonzero finitely generated fractional ideal $I$ of $D$ is $\star$--invertible
(i.e. ($II^{-1})^\star=D^\star$). In the semistar case, if $D$ is a P$\star$MD, then
the semistar integral closure of $D$ is integrally closed, thus, in the Krull's setting
of e.a.b. star operations, we recover the classical situation that $D$ has to be
integrally closed \cite{Houston/Malik/Mott:1984}. Several characterizations of
Pr{\"u}fer $\star$--multiplication domains were obtained recently, making also use of
the semistar Nagata  ring  $\Na(D,\star)$ and the semistar Kronecker functions ring
$\Kr(D,\star)$ (cf. \cite{fjs}, \cite{FL3}, \cite{BF} and \cite{Koch:2002}).

The starting point of this paper is the study of the localization of a P$\star$MD, $D$,
at any prime ideal (possibly, not quasi--$\star$--prime, i.e. not a prime ideal $P$
such that $P=P^\star\cap{D}$, since the localization at any prime ideal of this type is
known to be a valuation domain). As a consequence of this local study, we obtain new
examples of local P$\star$MDs.

One of the first results proved here is a characterization of a P$\star$MD, $D$,
through a local property, concerning the localizations of $D$ at a family of prime
ideals $P$ of $D$, and a global ``arithmetical" condition, concerning a finiteness
property of the ideals of the type $(aD\cap{bD})^\star$, see Theorem~\ref{cor:7}. We
apply this result to characterize P$\star$MDs as those domains such that the
localization at any prime ideal $P$ is a P$\star_P$MD (where $\star_P$ is a semistar
operation canonically associated to $\star$ by ``ascent" to $D_P$ ), see
Theorem~\ref{cor:3.7}. This result points out the important fact that Pr{\"u}fer
multiplication-like properties are really local properties and it opens the way for a
local--global study of Pr{\"u}fer $\star$--multiplication domains.

In order to realize these results we develop, preliminarily, a study on the behavior of
semistar operations properties under localizations. In particular, we show that finite
type, spectral, stable, a.b. or e.a.b. properties on a semistar operation $\star$
transfer to the induced semistar operation $\star_P$, defined on $D_P$, for any any
prime ideal $P\in\Spec(D)$.

At this stage, it is also natural to investigate on the relationship between the
semistar Nagata ring $\Na(D,\star)$ [respectively, semistar Kronecker function ring
$\Kr(D,\star)$] and the Nagata rings $\Na(D_P,\star_P)$ [respectively, the
Kronecker function rings $\Kr(D_P,\star_P)$] (i.e., on the local behavior of the
general Nagata ring and Kronecker function rings).  In this context, we
show that $\Na(D,\star)=\cap\{\Na(D_P,\star_P)\mid\;P\in\Spec(D)\}$ [respectively,
$\Kr(D,\star)=\cap\{\Kr(D_P,\star_P)\mid\;P\in\Spec(D)\}$]; we observe also that
the canonical inclusions $\Na(D,\star)_{D\setminus{P}}\subseteq\Na(D_P,\star_P)$
[respectively, $\Kr(D,\star)_{D\setminus{P}}\subseteq\Kr(D_P,\star_P)$] are not
equalities, in general.

In the last section we deal with the ``global" problem of building a semistar operation
in an integral domain $D$ by ``gluing" a given family of semistar operations
defined on $D_P$, for $P$ varying in a subset $\Theta$ of $\Spec(D)$. Since the
description of this semistar operation is in part folklore (at  least  in the star
setting), we deal specially with the problem of which properties, verified by all the
semistar operations defined on the localizations $D_P$, transfer to the ``glued''
semistar operation defined on $D$.  Among the other results, we prove that the finite
type and stable properties pass on, in the case the representation
$D=\cap\{D_P\;|\;P\in\Theta\}$ has finite character.  In order to glue semistar
operations verifying other relevant properties, like e.a.b. and a.b., we
evidentiate some obstructions; in fact, this type of semistar operations we prove
giving rise to a semistar operation of the same type under an extra condition, a sort
of ``stability under generalizations", denoted here by ($\downarrow$), see Theorem
\ref{pr:4.7}.  Finally, we have included several examples in order to better illustrate
the different constructions considered here and to show the essentiality of the
assumptions in the main results.

\section{Background}

Let us consider a commutative integral domain $D$ with quotient field $K$. Let
$\oF(D)$ [respectively, $\F(D)$ and $\f(D)$] denote the set of nonzero $D$--submodules
of $K$ [respectively, fractional ideals and nonzero finitely generated $D$--submodules
of $K$]. Note that $\f(D)\subseteq\F(D)\subseteq\oF(D)$.

\it A semistar operation on $D$ \rm is a map $\star\colon\oF(D)\longrightarrow\oF(D)$,
$E\mapsto{E^\star}$, such that
\bara
\item
$(xE)^\star=xE^\star$;
\item
$E\subseteq{F}$ implies $E^\star\subseteq{F^\star}$;
\item
$E\subseteq{E^\star}$ and $E^\star=(E^\star)^\star=: E^{\star\star}$,
\eara
for each $0\neq{x}\in{K}$ and for all $E$, $F\in\oF(D)$.

\noindent When $D^\star = D$, we say that $\star$ is \it a (semi)star operation on \rm
$D$. The \it identical operation $d_{D}$ on $D$ \rm (simply denoted by $d$), defined by
$E^{d_{D}} := E$, for each $E \in \oF(D)$, is a (semi)star operation on $D$. If not
stated explicitly, we generally assume that $\star$ is not the \it trivial semistar
operation $e_{D}$ on $D$ \rm (simply denoted by $e$), where $E^{e_{D}} := K$, for each
$E \in \oF(D)$. It is easy to see that $\star \neq e$ if and only if $D\neq{K}$ implies
that $D^\star\neq{K}$.

A semistar operation $\star$ is of \emph{finite type}, if for any $E\in\oF(D)$, we
have:
\[
 E^{\star_f}:=\cup\{F^\star\;|\;\,F\subseteq{E}\mbox{ and
}F\in\f(D)\}=E^\star.
\]
In general, for each semistar operation defined on $D$, $\star_f$, as defined before,
is also a semistar operation on $D$ and $\star_f\leq\star$, i.e., for any $E\in\oF(D)$
we have $E^{\star_f}\subseteq{E^\star}$.

There are several examples of finite type semistar operations; the most known is
probably the $t$--operation. Indeed, we start from the \it $v_{D}$--operation \rm on an
integral domain $D$ (simply denoted by $v$), defined as follows:
\[
 E^v :=(E^{-1})^{-1}=(D\colon(D\colon{E})) ,
\]
for any $E\in\oF(D)$, and we set $t_{D}:=(v_{D})_f$ (or, simply, $t =v_{f}$).

\noindent Other examples can be constructed as follows: let $T$ be an overring of an
integral domain $D$ and let $\star_{\{T\}}$ be the semistar operation on $D$, defined
by $E^{\star_{\{T\}}} := ET$, for each $E \in \oF(D)$, then $\star_{\{T\}}$ is a finite
type semistar operation on $D$. In particular, $\star_{\{D\}}= d$ and $\star_{\{K\}}=
e$.

In the following, we will construct several new semistar operations from a given one
and we will show that most of them are of finite type.

Let $\star$ be a semistar operation on $D$, we define the following set of prime ideals
of $D$:
\[
 \Pi^\star:=\{P\in\Spec(D)\;|\;\,P\neq (0) \mbox{ and
}P^\star\cap{D}\neq{D}\}.
\]
It may be that $\Pi^\star$ is an empty set but, in the particular case in which $\star$
is nontrivial of finite type, then $\Pi^\star \neq \emptyset$. In fact, we have that
each proper ideal $I$ is always contained in a proper \emph{quasi--$\star$--ideal} of
$D$, i.e. a proper ideal $J$ of $D$ such that such that $J^\star\cap{D}=J$, and
moreover (in the finite type case) each proper quasi--$\star$--ideal of $D$ is
contained in a maximal element in \it the set $\Gamma$ of all proper
quasi--$\star$--ideals; \rm finally, maximal elements in $\Gamma$ are prime ideals. We
will denote by
$$
 \calQ(\star):= \Spec^{\star}(D) :=\{Q \in \Spec(D)\;| \;\, Q =Q^{\star}\cap
D \}
$$
the set of all the \emph{ quasi--$\star$--prime ideals} of $D$.

If $\star$ is possibly not of finite type, we say that \it $\star$ possesses enough
primes\rm , if each proper quasi--$\star$--ideal is contained in a
quasi--$\star$--prime ideal. As a consequence, for each semistar operation that
possesses enough primes (e.g. a nontrivial finite type semistar operation), we have
$\Max(\Pi^\star)=\Max(\Gamma)$ (and it is a nonempty set). We will denote simply by
$\calM(\star)$ the set of all the maximal elements in $\Gamma$.

After analyzing this situation, we ask about how to relate semistar operations and
prime ideals. For any nonempty set $\Pi$ of prime ideals on $D$, we define a semistar
operation $\star_\Pi$ on $D$ as follows:
\[
 E^{\star_\Pi}:=\cap\{ED_P\;|\;\,P\in\Pi\},
\]
for any $E\in\oF(D)$.  If $\Pi=\emptyset$, we set $\star_{{\emptyset}}:= e$\ \hskip
-3pt;\  obviously, if $\Pi=\{(0)\}$, then $\star_{\Pi}:= e$\ \hskip -3pt. Semistar
operations defined by sets of prime ideals are called \emph{spectral semistar
operations}.

If $\star$ is a finite type semistar operation, we define
$\star_{sp}:=\star_{\calM(\star)}$; we have always that $\star_{sp}\leq\star$. Of
course, we can start from a general semistar operation $\star$, in that case we have a
new spectral semistar operation $\widetilde{\star}:=(\star_f)_{sp}$; this is the
biggest finite type spectral semistar operation in the set of all the finite type
spectral semistar operations on $D$, smaller or equal to $\star$.

\noindent A semistar operation $\star$ is \emph{stable} if
$(E\cap{F})^\star=E^\star\cap{F^\star}$, for each pair $E$, $F\in\oF(D)$; any spectral
semistar operation is stable, and stable finite type semistar operations coincide with
spectral finite type semistar operations.

The semistar operation $\widetilde{\star}$ can be also defined using the general
\emph{Nagata ring $\Na(D,\star)$ associated to} $\star$. Indeed, if we define:
\[
 \Na(D,\star):=D[X]_{N(\star)},
\]
where $N(\star):=\{f\in{D[X]}\;|\;\,\content(f)^\star=D^\star\}$ is a saturated
multiplicative subset of $D[X]$, then, for any $E\in\oF(D)$, we have
$E^{\widetilde{\star}}=E\Na(D,\star)\cap{K}$.

Nagata ring has a parallel behavior to the general \emph{Kronecker function ring}
associated to a semistar operation $\star$, defined as follows:
\[
\begin{array}{l}
 \Kr(D,\star)
 =\left\{\frac{f}{g}\in{K(X)}\;|\;\, f,\ g \in D[X] \setminus
 \{0\}\ \mbox{ and there exists } 0\neq{h}\in{D[X]}\right.\\
 \hspace*{22ex}\mbox{such that }
(\content(f)\content(h))^\star\subseteq(\content(g)\content(h))^\star\mbox{\LARGE{\}}}\cup\{0\}.

\end{array}
\]
Note that $\Na(D,\star)\subseteq\Kr(D,\star)$ and a ``new'' finite type semistar
operation can be defined on $D$ by the Kronecker function ring by setting:
\[
 F^{\star_a}:=F\Kr(D,\star)\cap{K},
\]
for any $F\in\f(D)$. This finite type semistar operation $\star_a$ on $D$ has another
more arithmetic description, as follows:
\[
 F^{\star_a}=\cup\{((FH)^\star\colon{H^\star})\;|\;\,H\in\f(D)\},
\]
for any $F\in\f(D)$. Moreover, $\star_a$ has a useful ``cancellation'' property: if
$E$, $F$, $G\in\f(D)$ and $(EF)^\star\subseteq(EG)^\star$, then
$F^\star\subseteq{G^\star}$. A semistar operation satisfying this property is called an
\emph{e.a.b. (= endlich arithmetisch brauchbar) semistar operation}. In the previous
``cancellation'' property, if we take $E\in\f(D)$ and $F$, $G\in\oF(D)$, a semistar
operation having this modified cancellation  property  is called an \emph{a.b. semistar
operation}. In general, we have the following characterizations: a finite type semistar
operation $\star$ is e.a.b. if and only if it is a.b. if and only if $\star=\star_a$.

A Kronecker function ring, and its ``counterpart'' the associated finite type a.b.
semistar operation, parameterizes certain valuation overrings of $D$. A valuation
overring $V\supseteq{D}$ is called a \emph{$\star$--valuation overring} of $D$ if
$F^\star\subseteq{FV}$ for any $F\in\f(D)$ (or, equivalently, if $\star_{f}\leq
\star_{\{V\}}$), then a finite type a.b. semistar operation $\star$ is characterized by
the following property:
\[
 F^\star=\cap\{FV\;|\;\,V\mbox{ is a $\star$--valuation overring of }D\},
\]
for each $F\in\f(D)$.

We refer to \cite{fjs}, and to the references contained in that paper, as a documented
source on semistar operations and on some of their properties briefly recalled above.

The following notation shall be use throughout the text. For a nonempty subset
$\Delta\subseteq\Spec(D)$ of prime ideals of an integral domain $D$ we define
\[
 {\Delta}^{\downarrow}:=
 \{Q\in\Spec(D)\;|\;\,Q\subseteq{P}\,,\mbox{ for some }P\in{\Delta}\},
\]
and we say that ${\Delta}$ is \emph{closed under generizations} if
${\Delta}={\Delta}^{\downarrow}$. In the same way, for any prime ideal $P\in\Spec(D)$,
we define $P^{\downarrow}:=\{P\}^{\downarrow}$.

\section{Localizing semistar operations}

Let $\star$ be a semistar operation on an integral domain $D$ and let $K$ be the
quotient field of $D$. For each $P \in \Spec(D)$, we consider the inclusion $D
\subseteq D_{P}$ of $D$ into its localization $D_{P}$ and the semistar operation $\sdDP
$, denoted simply $\sP$, on $D_{P}$, obtained from $\star$ by ``ascent to'' $D_{P}$,
i.e. $ E^{\sdDP} := E^{\sP} := E^\star$, for each $E \in
\boldsymbol{\overline{F}}(D_{P}) \, (\,\subseteq \boldsymbol{\overline{F}}(D)\,)$. Note
that if $P = (0)$ then $D_{P} = K$ and so $\sP$ coincides with $d_{K }\ (=e_{K})$ on
$K$.

\medskip
Our first goal is to study the transfer of some relevant properties from $\star$ to
$\star_P$.

\begin{lemma}\label{le:1}
Let $\star$ be a semistar operation on an integral domain $D$ and let $P\in\Spec(D)$.
\balf \item \sl If $\star$ is a finite type semistar operation on $D$, then $\sP$ is a
finite type semistar operation on $D_{P}$.  \item \sl If $\star$ is an e.a.b.
[respectively, a.b.] semistar operation on $D$, then $\sP$ is an e.a.b. [respectively,
a.b.] semistar operation on $D_{P}$.  \ealf
\end{lemma}
\begin{proof}
(a) is a consequence of \cite[Example 1 (e.1)]{fjs}.

(b) We give the proof in the a.b. case; a similar argument shows the e.a.b. case. Let
$G,\, H \in \boldsymbol{\overline{F}}(D_P)\subseteq \boldsymbol{\overline{F}}(D)$ and
$F \in \boldsymbol{f}(D_P)$ such that $(FG)^{\sP}\subseteq (FH)^{\sP}$. Since we can
find $F_0\in \boldsymbol{f}(D)$ such that $F=F_0D_P$, then we obtain
$(F_0G)^\star=(F_0D_PG)^{\sP}\subseteq(F_0D_PH)^{\sP}=(F_0H)^\star$. Therefore,
$G^{\sP} = G^\star \subseteq {H}^\star={H}^{\sP}$, because $\star$ is a.b.
\end{proof}

If $\star$ is a finite type a.b. semistar operation on an integral domain $D$, then it
is wellknown that $\star$ coincides with the semistar operation $\star_{\calV}$, where
$\calV:=\calV(\star):=\{V\supseteq{D}\;|\;\,V\mbox{ is a $\star$--valuation
overring}\}$ and $E^{\star_{\calV}}:=\cap\{EV\,|\;\,V\in\calV\}$, for each $E\in\oF(D)$
\cite[Lemma 2.8 (d)]{fjs}.

 \begin{corollary}\label{co:2.2}
  Let $\star$ be a finite type a.b. semistar operation  on an
 integral domain $D$.   For each prime ideal $P\in\Spec(D)$, we
 consider the (finite type a.b.) semistar operation $\star_{P}$ on
 $D_{P}$ and the subset of overrings
 $\calV_P:=\{V\supseteq{D}\;|\;\,V\mbox{ is a $\star$--valuation}$
 overring  $ \mbox{and } D_P\subseteq{V}\}$.  Then $\calV_P$ is exactly
 the set $\calV(\star_P)$ of the $\star_P$--valuation overrings of
$D_P$, i.e. $E^{\star_{P}}:=\cap\{EV\,|\;\,V\in\calV(\star_P)\}=\cap \{EV\,|\;\,V\in
\calV_P\}$, for each $E\in\oF(D_{P})$.
 \end{corollary}
 \begin{proof} We know already, from the previous lemma, that $\star_{P}$ is
 a finite type a.b. semistar operation on $D_{P}$. If $V\in\calV_P$ and
 $F\in\f(D)$, then $(FD_P)^{\star_P}=(FD_P)^\star\subseteq{(FD_P)V}$,
 thus $V$ is a $\star_P$--valuation overring of $D_{P}$. Conversely, if
 $W$ is a $\star_P$--valuation overring of $D_{P}$, then we have
 $F^\star\subseteq(FD_P)^\star=(FD_P)^{\star_P}\subseteq{D_PW}=FW$, for each
 $F\in\f(D)$, thus $W\ (\supseteq D_{P})$ is a $\star$--valuation
 overring of $D$. As a consequence, for each $F\in\f(D)$, we obtain that
 $(FD_P)^\star=(FD_P)^{\star_P}=\cap\{FD_PW\;|\;\,W \mbox{ is}$ $\mbox{ a
 $\star_P$--valuation
 overring}\}=\cap\{FV\;|\;\,V\in\calV_P\}=F^{\star_{\calV_P}}$. The
 conclusion follows since $\star_P$ is a finite type semistar operation.
 \end{proof}

\begin{proposition}\label{pr:2}
 Let $\star$ be a semistar operation on an integral domain $D$ and
let $P \in
 \Spec(D)$. If $\star = \tilde{\star}$ (is a finite type stable
 semistar operation on $D$), then $\sP =\widetilde{\sP}$ (is a finite
 type stable semistar operation on $D_{P}$).
\end{proposition}
\begin{proof}
Note that if $\star$ is a stable semistar operation then, from the definitions of
stability and of the semistar operation $\star_{P}$, it follows that $\star_{P}$ is a
stable semistar operation.  The conclusion follows from Lemma \ref{le:1} (a).

We give another proof that describes explicitly the set $\calQ(\star_{P})$ of all the
quasi--$\star_{P}$--prime ideals of $D_{P}$ in relation with the set $\calQ(\star)$ of
all the quasi--$\star$--prime ideals of $D$.

To avoid the trivial case, we can assume that $\star \neq e_{D}$ and that $P$ is a
nonzero prime ideal of $D$. Note that $\star =\star_{f}$, because $\star$ is a finite
type semistar operation, and so $\star_{P} =(\star_{P})_{f}$ (Lemma \ref{le:1} (a)).
Let $\calQ: = \calQ(\star)$ be the set of all the quasi--$\star$--prime ideals of $D$
then, for each $F \in \boldsymbol{\overline{F}}(D)$, we have $F^\star= \cap \{FD_Q \;|
\;\,Q\in\calQ\}$, because $\star = \tilde{\star}$ \cite[Corollary 2.11 (2)]{FL3}.

Assume that $P\in \calQ$. Let
$E\in\boldsymbol{\overline{F}}(D_P)\subseteq\boldsymbol{\overline{F}}(D)$,
 then
\[
 E^{\sP}=E^\star=\cap\{ED_Q\;| \;\,Q\in\calQ\}=
 \cap\{ED_Q\;| \;\,Q\in\calQ\setminus\{P\}\}\cap{E}\subseteq{E},
\]
hence $\sP=d_{D_{P}}$ is the identical (semi)star operation on $D_{P}$ and so,
obviously, $\sP=d_{D_{P}}= \widetilde{d_{D_{P}}} =\widetilde{\sP} $.

Assume that $P\notin\calQ$. Let
$E\in\boldsymbol{\overline{F}}(D_P)\subseteq\boldsymbol{\overline{F}}(D)$, then
\[
\begin{array}{{ll}}
 E^{\sP}&=(ED_P)^{\sP}=(ED_P)^\star=\cap\{ED_PD_Q\;| \;\,Q\in\calQ\}=\\
 &= (\cap\{ED_PD_Q\;| \;\,Q\in\calQ_{0}\})\cap(\cap\{ED_PD_Q\;|
\;\,Q\in\calQ_{1}\})
\end{array}
\]
where $\calQ_{1}:= \{Q \in \calQ\;| \;\, P\cap Q \mbox{ contains a nonzero prime ideal
of } D \}$ and $\calQ_{0}:= \{Q \in \calQ\;| \;\, P\cap Q \mbox{ does not contain a
nonzero prime ideal of } D \}$.

Note that if $Q\in \calQ_{0}$, i.e. if $P\cap Q$ does not contain a nonzero prime
ideal, then $D_PD_Q$ coincides necessarily with $K$, the quotient field of $D$.

Assume that $Q\in \calQ_{1}$. It is wellknown that there exists a natural bijective
correspondence between the set of prime ideals of $D_PD_Q$ and the set
$\{H\in\Spec(D)\;| \;\, H\subseteq{P}\cap{Q}\}$, hence $D_PD_Q=\cap\{D_H\;| \;\,
H\subseteq{P}\cap{Q}\mbox{ and }H\in\Spec(D)\}$. Moreover, note that the set $\calS(P;
Q)$ of all nonzero quasi--$\star$--ideals $I$ of $D$ contained in ${P}\cap{Q}$ is not
empty (since at least $H^{\star} \cap D$ is in $\calS(P; Q)$, if $H \subseteq
{P}\cap{Q} \mbox{ and } H\in \Spec(D)$). It is easy to see that the set of maximal
elements $\calS(P; Q)_{\mbox{\tiny max}}$ of $\calS(P; Q)$ is a nonempty set of prime
ideals, with $\calS(P; Q)_{\mbox{\tiny max}}\subseteq \calQ_{1} \subseteq\calQ$ and,
furthermore, each prime ideal $H$, with $H\subseteq {P} \cap {Q} $, is contained in
some element of $\calS(P; Q)_{\mbox{\tiny max}}$. Thus, we can rewrite:
\[
\begin{array}{{rl}}
 E^{\sP}\hskip -10pt
 &=EK \cap (\cap\{ED_H| \;H\subseteq{P}\cap{Q}\mbox{ and\!
}H\in\Spec(D),\mbox{ for $Q$
 varying in } \calQ_{1} \})\hskip -2pt=\\
 \hskip -10pt &=\cap\{ED_H\,|\; H\in \calS(P; Q)_{\mbox{\tiny max}}\,, \mbox{
for $Q$ varying in }
 \calQ_{1} \}\hskip -2pt=\\
\hskip -10pt &=\cap\{ED_H\,| \; H\subseteq{P} \mbox{ and } H\in\calQ\}.
 \end{array}
\]
\noindent Therefore it is easy to see that the set $ \{HD_P\;| \;\, H\subseteq{P}
\mbox{ and } H\in\calQ\}$ coincides with the set $\calQ(\star_{P})$ of all the
quasi--$\star_{P}$--prime ideals of $D_{P}$, which ``defines'' $\widetilde{\sP}$, i.e.
$\sP=\widetilde{(\sP)}$.
\end{proof}

\begin{remark}\label{rk:3}
\rm%
Note that the proof of the previous proposition shows the following statement:\
\sl%
If $\star$ is a finite type spectral semistar operation on $D$, defined by a subset
$\Delta \subseteq\Spec(D)$ (i.e. $\star:=\star_{{\Delta}}$), then $\sP$ is also a
finite type spectral semistar operation on $D_{P}$ and it is defined by the set
$\Delta_{P}:= \{HD_P\;| \;\,H\subseteq{P},\;H\in\Delta\}$ (i.e. $\star_{P
}=\star_{{\Delta_{P}}}$).
\rm%
\end{remark}

\begin{remark}\label{rk:31}
Let $\star$ be a semistar operation on $D$ and $P\in\Spec(D)$, then we have the
following diagram of semistar operations on $D_P$.
\[
 \begin{diagram}
 \star_P &\\
 \dLine &\\
 (\star_P)_f &\\
 \dLine &\rdLine &\\
 & &\widetilde{\star_P}\\
 & &\dLine\\
 ((\star_f)_P)_f=(\star_f)_P & &\\
 &\rdLine &\\
 & &\widetilde{(\star)}_P=\widetilde{(\widetilde{(\star)}_P)}
\end{diagram}
 \]
where the equalities are direct consequences of Lemma \ref{le:1} and Proposition
\ref{pr:2}. As a consequence of the proof of Proposition \ref{pr:2}, we obtain that:
$$
\begin{array}{rl}
 P\in\Spec^{\star_f}(D)\,\;\; \Rightarrow\,\,
 &\Spec^{\widetilde{\star_P}}(D_P) = \Spec^{(\star_P)_f}(D_P)=\\
 & = \Spec^{(\star_f)_P}(D_P) =\Spec^{\widetilde{(\star)}_P}(D_P)\,.
\end{array}
$$
Indeed, if $QD_P\in\Spec^{\widetilde{(\star)}_P}(D_P)$, then
\[
\begin{array}{llll}
 QD_P
 &=(QD_P)^{\widetilde{(\star)}_P}\cap{D_P}=
(QD_P)^{\widetilde{\star}}\cap{D_P}=\\
 &=(\cap\{QD_PD_H \mid\, H\in\Spec^{\star_f}(D)\}) \cap D_P =
 (QD_P)^{\widetilde{\star_P}}\cap D_P \subseteq\\
 &\subseteq{QD_P}\cap{D_P}=QD_P\,.
\end{array}
\]
\end{remark}

\medskip

Let $\star$ be a semistar operation on an integral domain $D$. Assume that $D$ is a
P$\star$MD, i.e. an integral domain such that each $F\in \boldsymbol{f}(D)$ is a
$\star_f$--invertible ideal, that is ${(FF^{-1})}^{\star_f} = D^{\star}$. P$\star$MDs
are characterized in several different ways in \cite[Theorem 1]{fjs}: for instance,
$D_Q$ is a valuation domain for each $Q\in\calM(\star_f)$, where $\calM(\star_f)$ is
the (nonempty) set of all the maximal elements of $\Spec^{\star_f}(D)$.  A consequence
of this fact is that $D$ is a P$\star$MD if and only if it is a P$\tilde{\star}$MD,
since $\calM(\tilde{\star})=\calM(\star_f)$ \cite[Lemma 3 (g)]{fjs}.

\begin{remark} \label{rk:4}
Note that, \sl if $0\neq{a}$, $b$ belong to an integral domain $D$, then the following
equality holds:
$$
\frac{1}{ab}(aD\cap{bD})=((a,b)D)^{-1}.
$$
 \rm
\end{remark}
\begin{proof}
Let $x\in{aD}\cap{bD}$, then $x=ax'=bx''$, for some $x'$, $x''\in{D}$, thus we obtain
that $\frac{x}{ab}=\frac{1}{b}x'=\frac{1}{a}x''$. Henceforth, for each
$az_1+bz_2\in(a,b)D$, we have $\frac{x}{ab}(az_1+bz_2)=x''z_1+x'z_2\in{D}$. Therefore
$\frac{x}{ab}\in ((a,b)D)^{-1}$.\ Conversely, let $y\in((a,b)D)^{-1}$, then $ya=x'$ and
$yb=x''$, for some $x'$, $x''\in{D}$. Henceforth $aby=bx'=ax''$, and so
$aby\in(aD\cap{bD})$.
\end{proof}

\begin{lemma} \label{le:5}
Let $\star$ be a semistar operation on an integral domain $D$. Assume that $D$ is a
P$\star$MD. For each pair of nonzero elements ${a}$, $b\in{D}$, there exists $F\in
\boldsymbol{{f}}(D)$ such that:
$$
(aD\cap{bD})^\star = F^\star
$$
(\,in this situation, we say briefly that $(aD\cap{bD})^\star$ is an\ \rm ideal of
$\star$--finite type\,\sl ).
\end{lemma}
\begin{proof}
Recall that $\frac{1}{ab}(aD\cap{bD})=((a,b)D)^{-1}$ and thus note that
$(aD\cap{bD})^\star$ is an ideal of $\star$--finite type if and only if
$(((a,b)D)^{-1})^\star$ is \it a (fractional) ideal of $\star$--finite type, \rm i.e.
there exists $G\in \boldsymbol{{f}}(D)$ such that $G^\star = (((a,b)D)^{-1})^\star$,
(Remark \ref{rk:4}). Since $D$ is a P$\star$MD, each nonzero finitely generated
(fractional) ideal $F$ of $D$ is $\star_{f}$--invertible and so $(F^{-1})^\star$ is
also a (fractional) ideal of $\star$--finite type. Therefore, we conclude that
$(((a,b)D)^{-1})^\star$ is a (fractional) ideal of $\star$--finite type.
\end{proof}

The converse of this result holds if we add some extra conditions.

\begin{theorem} \label{th:6}
Let $\star$ be a spectral semistar operation on $D$ defined by a set $\Delta$ of
valuation prime ideals of $D$ (or, equivalently, by a family of essential valuation
overrings of $D$), i.e. $E^\star:=\cap\{ED_P\;|\,\;P \in \Delta \}$, for each $E\in
\boldsymbol{\overline{F}}(D)$, where $D_{P}$ is a valuation domain for each $P \in
\Delta$. If, for each pair of nonzero elements $a$, $b\in{D}$, we have that
$(aD\cap{bD})^\star$ is an ideal of $\star$--finite type, then $D$ is a P$\star$MD.
\end{theorem}
\begin{proof}
{Step 1}. \sl If $I$ is a finitely generated ideal of $D$ and $a\in{D}$, with $a \neq
0$, then $(I\cap{aD})^\star$ is an ideal of $\star$--finite type. \rm

The proof is based on an argument from \cite{Zafrullah:fc}. Set $I:=\sum_{i=1}^na_iD$
and $J_i:=a_iD\cap{aD}$, for $1 \leq i \leq n$. By the hypothesis $(J_i)^\star$ is an
ideal of $\star$--finite type. For each index $i$, $1 \leq i \leq n$, let $F_i$ be a
finitely generated ideal of $D$ such that $(J_i)^\star=(F_i)^\star$. Since $D_P$ is a
valuation domain and a $D$--flat overring of $D$, we have:
\[
\begin{array}{ll}
 (\sum_i F_i)^\star &={(\sum_i {(F_i)}^\star)}^\star= (\sum_i
 (J_i)^\star)^\star = (\sum_iJ_i)^\star =\\
 &=\cap_{ P\in \Delta }\left(\sum_i (a_iD \cap {aD})\right)D_P =\cap_{P\in
 \Delta }\left(\sum_i(a_iD\cap{aD})D_P\right) = \\
 &=\cap_{P\in \Delta }\left(\sum_i(a_iD_P\cap{aD_P})\right)
 =\cap_{P\in \Delta }\left((\sum_ia_i)D_P \cap{aD_P}\right)= \\
 &=\left(\cap_{P\in \Delta }(\sum_ia_i)D_P\right) \cap \left(\cap_{P\in
\Delta }aD_P\right)
 =(\sum_ia_iD)^\star \cap (aD)^\star =\\
 &=I^\star \cap (aD)^\star =(I\cap{aD})^\star \qquad \qquad\mbox{(as $\star$
is stable).}
\end{array}
\]

{Step 2}. \sl Any finite intersection of nonzero principal ideals of $D$ is an ideal of
$\star$--finite type. \rm

Let $a_1, a_2,\ldots, a_t \in D$ be a family of nonzero elements. We may assume that $t
\geq 2$ and, by induction on $t$, we may assume that $a_1D \cap a_2D \cap \ldots \cap
{a_{t-1}D}$ is an ideal of $\star$--finite type, i.e. there is a finitely generated
ideal $F$ of $D$ such that $(a_1D \cap a_2D\cap\ldots\cap{a_{t-1}D})^\star=F^\star$.
Then, we have:
\[
\begin{array}{rl}
 (a_1D\cap a_2D \cap\ \ldots\ \cap{a_{t}D})^\star =& \left((a_1D
\cap a_2D \cap\ \ldots\ \cap{a_{t-1}D})^\star \cap(a_{t}D)^\star
\right)^\star =\\
=& \left(F^\star
 \cap (a_{t}D)^\star\right)^\star =\left(F\cap a_{t}D \right)^\star
 \end{array}
\]
and, by Step 1, this is an ideal of $\star$--finite type.
\smallskip

{Step 3}. \sl If $I$ is a nonzero finitely generated ideal of $D$, then $I^{-1}$ is a
(fractional) ideal of $\star$--finite type. \rm

 The case of an ideal $I$ generated by two elements $ a, b \in D$ follows
 immediately from Remark~\ref{rk:4}, since we know already that
 $(aD\cap{bD})^\star$ is an ideal of $\star$--finite type if and only if
 $(((a,b)D)^{-1})^\star$ is a (fractional) ideal of $\star$--finite type. The
 conclusion follows from the assumption that
 $(aD\cap{bD})^\star$ is an ideal of $\star$--finite type. \;
 The general case of a finitely generated ideal $I := (x_1, x_2, \ldots,
 x_{t})D$ follows from Step 2. In fact, without loss of generality,
 we can assume that $x_{i}\neq 0$, for each $1 \leq i \leq t$, thus:
 $$
 I^{-1} = (D: I) = (D : \sum_{1\leq i \leq t\ }x_{i}D) = \cap_{1\leq i \leq
 t\ } (D:x_{i}D) = \cap_{1\leq i \leq t\ } x^{-1}_{i}D,
 $$
 and, if we write $x^{-1}_{i}:= a_{i}/d$, with $a_{i}$ and $d$ nonzero
elements
 in $D$, for $1 \leq i \leq t$, then:
 $$
 (\cap_{1\leq i\leq t\ } x^{-1}_{i}D)^\star \hskip -2pt =(d^{-1}(\cap_{1\leq i
 \leq t\ } a_{i}D))^\star \hskip -2pt  = d^{-1}(\cap_{1\leq i \leq t\ }a_{i}D)^\star
 \hskip -2pt = d^{-1}F^\star \hskip -2pt = (d^{-1}F)^\star,
 $$
 for some $F\in\f(D)$.

Let $I$ be any nonzero finitely generated ideal of $D$. By Step 3, we know that
$I^{-1}$ is a (fractional) ideal of $\star$--finite type. Since $D_{P}$ is a valuation
domain and a $D$--flat overring of $D$, we have:
\[
\begin{array}{ll}
 (I\;I^{-1})^\star
 &=\cap_{P\in \Delta }(I\;I^{-1})D_P = \cap_{P\in \Delta
}(ID_P\;I^{-1}D_P)=\\
 &=\cap_{P\in \Delta }(ID_P\;(ID_P)^{-1})=\cap_{P\in \Delta }D_P =D^\star.
\end{array}
\]
We conclude that $D$ is a P$\star$MD.
\end{proof}

\begin{theorem} \label{cor:7}
Let $\star$ be a semistar operation on an integral domain $D$. Then the following
statements are equivalent:
\brom
\item \sl
$D$ is a P$\star$MD;
\item \sl
the following two conditions hold:
 {\bara
 \bf\item \sl
 for each subset $\Theta \subseteq\Spec(D)$, such that $\tilde{\star }=
 {{\wedge}}_{\widetilde{\Theta}}:=
 \wedge \{ \widetilde{\mbox{\d{$\star$}}_{P}}\,|\; P\in \Theta \}$,
 where $E^{{{\wedge}}_{\widetilde{\Theta}}} := \cap \{
(ED_{P})^{\widetilde{\star_{P}}} \,|\; P \in \Theta \}$, for each $E\in
\boldsymbol{\overline{F}}(D)$, we have that $D_P$ is a P$\sP$MD, for each $P \in
\Theta$;
 \bf\item \sl
 for any pair of nonzero elements $a$, $b\in{D}$, we have that
 $(aD\cap{bD})^\star$ is an ideal of $\star$--finite type.
 \eara}
 \erom
\end{theorem}
\begin{proof}
(i) $\Rightarrow$ (ii).  We can assume that $\star = \tilde{\star}$, since the notions
of P$\star$MD and P$\tilde{\star}$MD coincide \cite[Section 3, Theorem 1]{fjs} and
$\star_{P }= (\widetilde{\star})_{P}= \widetilde{\star_{P}}$ (Proposition~\ref{pr:2}).
By Lemma~\ref{le:5}, we only need to show that condition (a) holds.  More generally, we
show that, under the assumption (i), $D_P$ is a $P\sP$MD, for each $P\in\Spec(D)$. Let
$ID_P$ be a finitely generated ideal of $D_P$, with $I$ a nonzero finitely generated
ideal of $D$. Since $D$ is a P$\star$MD by assumption, there exists a finitely
generated (fractional) ideal $J$ of $D$ such that
$(IJ)^{\star_{f}}=(IJ)^\star=D^\star$.

Assume that $P \in \calQ({\star_{f}})=\calQ({\star})$\ (since $\star = \tilde{\star}$
implies that $\star = \star_{f}$). Then $E^\star D_P=ED_P$, for every $E\in
\boldsymbol{\overline{F}}(D)$ (because $\star = \widetilde{\star}$) hence, in
particular, $(IJ)^\star D_{P}= IJD_P=ID_PJD_P= D_P$ and thus $ID_P$ is invertible in
$D_P$ (i.e. $D_{P}$ is a valuation domain and so, trivially, it is a P$\ast$MD, for
every semistar operation $\ast$ on $D_{P}$).

Assume that $P\in \Spec(D) \setminus\calQ({\star})$. Then $(IJ)^\star=D^\star$ implies
$1\in{IJD_Q}$, for each $Q\in\calQ({\star})$; in particular $1\in{IJD_Q},$ for each
$Q\in\calQ({\star})$ such that $Q \subseteq {P}$ (this set of prime ideals is nonempty,
since each minimal prime ideal of a nonzero principal ideal of $D$ is in
$\calQ({\star})$, for any finite type semistar operation $\star$). Therefore, by the
proof of Proposition \ref{pr:2},
$D\subseteq\cap\{{IJD_Q}\mid{Q}\in\calQ({\star})\,,\;Q\subseteq{P}\}=
{(IJD_P)}^{\widetilde{\sP}}\subseteq(ID_P JD_P)^{\sP}= (ID_PJD_P)^{({\sP})_{f}}$, thus
we obtain $(ID_PJD_P)^{({\sP})_{f}}=D^{\sP}$ and, hence, $ID_P$ is a
$(\sP)_f$--invertible ideal of $D_{P}$ (with $\sP= (\sP)_f$), i.e. $D_P$ is a P$\sP$MD.
See also \cite[ Section 3, Theorem 1]{fjs}.

(ii) $\Rightarrow$ (i). Note that, by (a), we have that, for each $E\in
\boldsymbol{\overline{F}}(D)$, $E^{\tilde{\star}} = E^{{{\wedge}}_{\widetilde{\Theta}}}
= \cap \{ (ED_{P})^{\widetilde{\star_{P}}} \,|\; P \in \Theta \}= \cap \mbox{\Large \{}
\, \cap \{ (ED_{P})D_{Q}= ED_{Q} \,|\; Q \in {\calM}(\widetilde{\star_{P}}) \} \,\, |
\; P \in \Theta \mbox{\Large \}}$, where $D_{Q}$ is a valuation domain, for each $Q \in
\calM(\widetilde{\star_{P}})$ and for each $P \in \Theta$. By Theorem~\ref{th:6}, we
deduce that $D$ is a P$\tilde{\star}$MD, i.e. $D$ is a P${\star}$MD \cite[Section 3,
Theorem 1]{fjs}.
\end{proof}

\begin{remark}\label{rk:8}
 From the previous proof it follows that:\ \sl if $D$ is a P$\star$MD, then
 $D_{P}$ is a P$\sP$MD, for each $P \in \Spec(D)$. \rm In the next section,
 we will show that the converse holds. Furthermore, in Section 4, we will
 deepen the study of the semistar operations of the type
 $\wedge_{\Theta}$; in particular, we will establish a natural
 relation between the semistar operation ${\wedge}_{\widetilde{\Theta}}$
 (considered in Theorem~\ref{cor:7}) and the finite type stable
 semistar operation, $\widetilde{({\wedge}_{\Theta})}$, canonically
 associated to ${\wedge}_{\Theta}$, where ${\wedge}_{\Theta}$ is defined as
follows:
 $E^{{\wedge}_{\Theta}} := \cap \{ (ED_{P})^{\star_{P}} \,|\; P \in
 \Theta \}$, for each $E\in \boldsymbol{\overline{F}}(D)$.
 \end{remark}

\section{Compatibility with localizations}

Let $D$ be an integral domain with quotient field $K$ and let $P\in\Spec(D)$. On the
localization $D_{P}$ of $D$ at $P$, we can consider the (semi)star operation $v_{D_P}$
[respectively, the semistar operation $v_{P}:= {\dot{v}_{D}}^{D_{P}}$] which denotes
the (semi)star $v$--operation on $D_{P}$ [respectively, the semistar operation on $D_P$
induced by the the (semi)star $v$--operation $v_{D}$ on $D$]. If the conductor $(D :
{D_P})$ is zero, then $(D_P)^{v_P} = (D_P)^{v_{D}} = (D : (0)) = K$, hence, in general,
the (semi)star operation $v_{D_P}$ (on $D_{P}$) does not coincide with the semistar
operation $v_P$ (on $D_{P}$).

Let us now relate $v_P$ and $v_{D_P}$ in some particular case.

As a special case of \cite[Lemma 3.4 (2)]{Kang:1989} we have the following:

\begin{lemma} \label{le:3.1}
Let $D$ be an integral domain. For each $F\in\boldsymbol{f}(D)$ and for each $P \in
\Spec(D)$, we have $(FD_P)^{v_{D_P}}=(F^{v_D}D_P)^{v_{D_P}}$.
\end{lemma}
\vskip-.5cm \hfill $\square$

It is known that, if $D$ is a P$v$MD, then $D$ is a $v$--\emph{coherent} domain in the
sense of \cite{Fontana/Gabelli:1996}, i.e. if $I$, $J \in\boldsymbol{f}(D)$, then
$I^v\cap{J^v}$ is an ideal of $v$--finite type. Also in \cite{Fontana/Gabelli:1996}
there is the following characterization of $v$--coherent domains: $D$ is a
$v$--coherent domain if and only if, for each $I\in \boldsymbol{f}(D)$, there exists
$F\in \boldsymbol{f}(D)$ such that $I^{-1}=F^v$ (i.e. $I^{-1}$ is an ideal of
$v$--finite type).

\begin{lemma} \label{le:3.2}
Let $D$ be a $v$--coherent domain (in particular, a P$v$MD). For each
$F\in\boldsymbol{f}(D)$ and for each $P \in \Spec(D)$, we have
$(FD_P)^{v_{D_P}}=F^{v_D}D_P $.
\end{lemma}
\begin{proof} Since we are assuming that $D$ is a $v$--\emph{coherent}
domain, if $I\in \boldsymbol{f}(D)$, then there exists $F\in \boldsymbol{f}(D)$ such
that $I^{-1}=F^{v_D}$ (or, equivalently, $I^{v_D}=F^{-1}$). Now we localize both sides
of the previous equality at $P$ and we obtain:
\[
 F^{v_{D}}D_P= I^{-1}D_P=(ID_P)^{-1}.
\]
Since $I^{v_{D}} = (D:F)= (D:F)^{v_{D}}$, then $(FD_P)^{v_{D_{P}}}= (D_{P}:
(D_{P}:FD_P)) = (D_{P}: (D : F) D_P) = (D_{P}: I^{v_{D}} D_P) \subseteq (D_{P}: ID_P)=
(D: I)D_P = I^{-1}D_{P} = F^{v_{D}}D_P $. By the previous Lemma~\ref{le:3.1}, we know
that $(FD_P)^{v_{D_P}}=(F^{v_D}D_P)^{v_{D_P}}$, thus we conclude immediately that
 $(FD_P)^{v_{D_P}}=F^{v_D}D_P $.
\end{proof}

\begin{proposition} \label{pr:3.3}
If $D$ is a $v$--coherent domain (in particular, if $D$ is a P$v$MD domain) then, for
each $P \in \Spec(D)$ and for each $F\in\boldsymbol{f}(D_{P})$, we have
$F^{v_{D_P}}\subseteq F^{v_P} $.
\end{proposition}
\begin{proof}
For each $F\in\boldsymbol{f}(D_{P})\,(\subseteq \boldsymbol{\overline{F}}(D))$ there
exists $F_0\in \boldsymbol{f}(D)$ such that $F=F_0D_P$, then by using
Lemma~\ref{le:3.2} we have:
\[
\begin{array}{ll}
 F^{v_{D_P}} &= (F_0D_P)^{v_{D_P}} ={F_0}^{v_{D}}D_P \subseteq\\
 & \subseteq F^{v_{D}}D_P=F^{v_{D}}= (D: (D: F)) = F^{v_P}
 \end{array}
\]
(note that the first equality in the second line is a consequence of the following
general fact: if $E\in \boldsymbol{\overline{F}}(D_{P})$, then $(D : E)$ is also in
$\boldsymbol{\overline{F}}(D_{P})$ and so $E^{v_{D}}$ belongs to
$\boldsymbol{\overline{F}}(D_{P})$).
 \end{proof}

\begin{remark} \label{pb:3.4}
Let $\star$ be a semistar operation on $D$ and, for each $P\in\Spec(D)$, let $\sP$ be
the semistar operation induced on $D_P$, defined in the previous section.  For which
properties \textbf{(P)} concerning $(D, \star)$ we have that $(D_{P}, \sP)$ satisfies
\textbf{(P)} ?

 A positive answer to this question was already given for the following
 properties:\ \bf (a) \rm\ $\star$ is a finite type semistar operation on
 $D$;\ \bf (b) \rm\ $\star$ is a   stable semistar operation
 on $D$;\  \bf (c) \rm\ $\star$ is a finite type spectral
semistar
 operation on $D$;\ \bf (d) \rm\ $\star$ is a finite type stable
 semistar operation on $D$ (i.e. $\star=\tilde{\star}$);\ \bf (e) \rm\
 $\star$ is an e.a.b semistar operation on $D$; \bf (f) \rm\ $\star$
is
 an a.b semistar operation on $D$;\ \bf (g) \rm\ $D$ is a
 P$\star$MD, (cf.  Lemma \ref{le:1}, Proposition \ref{pr:2},
 Remark \ref{rk:3} and Remark \ref{rk:8}).
\end{remark}

  In this  ambit,  a natural problem is to study the behavior of the
 generalized Nagata ring and  of the generalized  Kronecker function
 ring in relation with the localization at any prime ideal $P$.   We have
 the following:

\begin{proposition} \label{pr:3.5} Let $\star$ be a semistar
operation on an integral domain $D$ and let $P \in \Spec(D)$. Then the following
statements hold:
\bara
\item\sl
$\Na(D,\star)_{D\setminus P}\subseteq\Na(D_{P},\star_{P})$.
\item\sl
$\Na(D,\star)=\cap\{\Na(D_{P},\sP)\,|\; P\in\Spec(D)\} =$ \\
\hspace*{9.5ex}$ =\cap\{\Na(D_{M},\sM)\,|\; M\in\Max(D)\}.$ \item\sl
$\Kr(D,\star)_{D\setminus{P}}\subseteq\Kr(D_{P},\star_{P})$. \item\sl
$\Kr(D,\star)=\cap\{\Kr(D_{P},\sP)\,|\; P\in\Spec(D)\}=$\\
 \hspace*{9.5ex}$=\cap\{\Kr(D_{M},\sM)\,|\; M\in\Max(D)\}.$
\eara
\end{proposition}

\begin{proof} (1) Set ${\calQ} := {\calQ}({\star_{f}})$. Recall from
the proof of Proposition~\ref{pr:2}, see Remark \ref{rk:31}, that
${\calQ}({(\star_{f})}_{P}) = \{ Q{D_{P}}\,|\; Q \subseteq P\,, \; {Q\in{\calQ}}\}$.
Set simply ${\calQ}_{P}:={\calQ}({(\star_{f})}_{P})$ and
$\mbox{\d{${\calQ}$}}{_P}:=\{Q\in\Spec(D)\,|\;QD_{P}\in{\calQ}_{P}\}$. Note also that $
\calQ = \cup \{ {\mbox{\d{${\calQ}$}}{_P}}\,|\; P \in \Spec(D) \}$. We know that
$\Na(D, \star) {= \Na(D, \star_{f}) }$ $= \cap \{D_{Q}(X)\,|\; Q \in {\calQ}\}$.
Therefore $\Na(D, \star)_{{D \setminus P}} = (\cap \{D_{Q}(X)\,|\; Q \in {\calQ} \})_{D
\setminus P} \subseteq \cap \{D_{Q}(X)_{D \setminus P}\,|\; Q \in {\calQ} \} = \cap
\{D_{Q}(X) \,|\;Q \in {\calQ}\ \mbox{ and } \, \ Q \subseteq P \} = \cap
\{D_{Q}(X)\,|\; $ $Q \in {\calQ}_{P}\} {= \Na(D_{P}, (\star_{f})_{P})}$ $ = \Na(D_{P},
\star_{P})$.

(2) Since $\calQ = \cup \{ {\mbox{\d{${\calQ}$}}{_P}}\,|\; P \in \Spec(D) \}$, then
$\Na(D, \star) = \cap \{D_{Q}(X)\,|\; Q \in {\calQ}\} = \cap \mbox{\Large{\{}} \, \cap
\{D_{Q}(X)\,|\; Q \in {\calQ_{P}}\}\ \,|\; P \in \Spec(D) \mbox{\Large{\}}} = \cap \{
\Na(D_{P}, \sP)\,|\; P \in \Spec(D)\}$. The proof is similar for the $\Max(D)$ case.

(3)  We start by recalling, from Corollary \ref{co:2.2}, the following fact:  \sl A \
$\sP$--\-valuation overring \ $W$ of $D_{P}$ is the same as a $\star$--valuation
overring $W$ of $D$ such that $W \supseteq D_{P}$.  \rm

We know that $\Kr(D, \star) = \cap \{W(X) \,|\; W \mbox{ is a $\star$--valuation }$
overring of  $D \}$. The\-refore, using  Corollary \ref{co:2.2}, the fact that
$W(X)_{D\setminus P}= W_{D\setminus P}(X)$ and that $W_{D\setminus P}$ is a valuation
overring of $D_{P}$, for each valuation overring $W$ of $D$, then
\[
\begin{array}{lll}
\Kr(D,\star)_{{D\setminus{P}}}
 &\hskip -6pt =(\cap\{W(X)\,|\;W\mbox{ is a $\star$--valuation overring of
}D\})_{D\setminus{P}}\subseteq\\
 &\hskip -6pt \subseteq\cap\{W(X)\,|\;W\mbox{ is a $\star$--valuation
 overring of }D\mbox{ and }W\supseteq{D_{P}}\}\!=\\
 &\hskip -6pt =\cap \{W(X)\,|\;W\mbox{ is a $\star_{P}$--valuation overring
of }D \}=\\
 &\hskip -6pt =\Kr(D_{P},\star_{P}).
\end{array}
\]

(4) From the Corollary \ref{co:2.2} we deduce that $\{W \,|\; W \mbox{ is a
$\star$--valuation}$ overring of   $ D \} = \cup \, \mbox{\Large{\{}} \, \{ W \, | \; W
\mbox{ is a $\sP$--valuation overring of } D_{P} \}\; | \; \, P \in \Spec(D)
\mbox{\Large{\}}}$. Therefore, we have that:
\[
\begin{array}{ll}
\Kr(D,\star)
 \hskip -9pt&=\cap\{W(X)\,|\;W\mbox{ is a $\star$--valuation overring
 of }D\}=\\
 \hskip -9pt&=\cap \mbox{\Large{\{}} \cap\!\{W(X)\,|\;W\mbox{ is a
 $\sP$--valuation overring of }D_{P} \}\,|\,
 P\in\Spec(D) \mbox{\Large{\}}}\\
 \hskip -9pt&=\cap\{\Kr(D_{P},\sP)\,|\;P\in\Spec(D)\}.
\end{array}
\]
The proof is similar for the $\Max(D)$ case.
\end{proof}

Next problem is to relate the Nagata ring or the Kronecker function ring, associated to
a localized semistar operation, to the corresponding localization of the Nagata ring or
of the Kronecker function ring, respectively. More precisely,

\begin{problem} \label{pb:3.6}\hfill
 \sl Let $D$ be an integral domain, $\star$ a semistar operation on
 $D$ and $P$ a prime ideal of $D$.
\bara
\item
 Under which conditions on $D$ and $P$, $\Na(D,\star)_{D\setminus
P}=\Na(D_{P},\star_{P})$ ?
\item
 Under which conditions on $D$ and $P$, $\Kr(D,\star)_{D\setminus
P}=\Kr(D_{P},\star_{P})$ ?
\item
 In case of a Pr{\"u}fer--$\star$--multiplication domain is the answer to
both questions positive ?
\eara
\end{problem}

We deal first with the Nagata ring. Without loss of generality, we may assume that
$\star=\widetilde{\star}$. For each prime ideal $P\in\Spec(D)$, we have the following
picture of finite type stable semistar operations:
\[
 \begin{diagram}
 D &\rTo & D[X] &\rTo^{\;\;\; \varphi_{_{P}}} && D_P[X]\\
 \star &\rMapsto & \eta\ (=N(\star)) &\rMapsto &&
 \mbox{$\left\{\begin{array}{l} \eta_P\ (=N(\star_P))\\
 \varphi_{_{P}}(\eta)\end{array}\right.$}
\end{diagram}
 \]
where $\eta$ [respectively, $\eta_P$] is the semistar operation on $D[X]$
[respectively, on $D_{P}[X]$] ``defined by the saturated multiplicative subset''
$N(\star)$ [respectively, $N(\star_P)$], i.e. $E^\eta := ED[X]_{N(\star)}= E\Na(D,
\star)$, for each $E \in \boldsymbol{\overline{F}}(D[X])$ [respectively, $E^{\eta_{P}}
:= ED_{P}[X]_{N(\star_{P})}= E\Na(D_{P}, \star_{P})$, for each $E \in
\boldsymbol{\overline{F}}(D_{P}[X])$], and $\varphi_{P}(\eta)$ (or, simply,
$\varphi(\eta)$) is the semistar operation on $D_{P}[X]$ defined as follows:
\[
 E^{\varphi(\eta)}:=\left\{z\in{K(X)}\;|\;\,\varphi^{-1}\left(E :_{D_{P}[X]\
 }{zD_{P}[X]}\right)\cap{N(\star)} \neq\emptyset\right\},
\]
for each $E \in \boldsymbol{\overline{F}}(D_{P}[X])$. Note that $\eta$ , $\eta_P$ and
$\varphi(\eta)$ are finite type stable semistar operations. In general, we have
$\varphi(\eta)\leq \eta_P$. Indeed, given an ideal $J\subseteq{D_P[X]}$ satisfying
$\varphi^{-1}(J)\cap{N(\star)}\neq \emptyset$, then there exists $f\in \varphi^{-1}(J)$
such that $\content(f)^\star=D^\star$, hence
$(\content(\varphi(f)))^{\star_{P}}=(\content(f)D_P)^{\star_{P}}
=(\content(f)D_P)^\star=(\content(f)^\star{D_P})^\star=(D^\star
D_P)^\star=D_P^\star=D_P^{\star_{P}}$, thus $\varphi(f) \in N(\star_{P})$. From this
remark, we deduce immediately that $E^{\varphi(\eta)} \subseteq E^{\eta_{P}}$, for each
$E \in \boldsymbol{\overline{F}}(D_{P}[X])$.

With this background, now we use some wellknown fact of hereditary torsion theories
\cite{BJV} or, equivalently, of localizing systems associated to semistar operations
\cite{FH}. More precisely, we know that applying a finite type stable semistar
operation is exactly the same as doing the localization with respect to the associated
finite type hereditary torsion theory or with respect to the associated localizing
system of ideals and, moreover, it is wellknown that it is possible to ``interchange'',
in a natural way, two subsequent localizations of the previous type. Therefore:
\[
\begin{array}{rl}
 \Na(D,\star)_{D\setminus{P}}
 &=\left(D[X]^\eta\right)_{D\setminus{P}}
 = \left(D[X]_{D\setminus{P}}\right)^{\varphi(\eta)}
 =\left(D_P[X]\right)^{\varphi(\eta)}\subseteq\\
 &\subseteq{D_P[X]^{\eta_P}}=\Na(D_P,\star_P),
\end{array}
\]
since the localizing system of ideals of $D[X]$ associated to $\eta$ is the set
${\calF}^{\eta}:= \{I \mbox{ ideal of } D[X] \mid \, I^\eta =D[X]^\eta \}= \{I \mbox{
ideal of } D[X] \mid \, I\cap N(\star) \neq \emptyset \}$ and the localizing system of
ideals of $D_P[X]$ associated to $\varphi(\eta)$ [respectively, $\eta_{P}$] is the set
${\calF}^{\varphi(\eta)}:= \{ID_{P}[X] \mid \, I \mbox{ ideal of } D[X],\
(ID_{P}[X])^{\varphi(\eta)} =D_{P}[X]^{\varphi(\eta)}\} = \{ID_P[X] \mid \, I \mbox{
ideal of } D[X],\ ID_P[X] \cap N(\star) \neq \emptyset \}$ [respectively,
${\calF}^{\eta_{P}}:= \{J \mbox{ ideal of } D_{P}[X] \mid \, J^{\eta_{P}}
=D_{P}[X]^{\eta_{P}}\}= \{J \mbox{ ideal of } D_{P}[X] \mid \, J\cap N(\star_{P}) \neq
\emptyset \}$] .

To show that $D_P[X]^{\varphi(\eta)}$ and ${D_P[X]^{\eta_P}}$ are either equals or
different we only need to compare the prime spectra associated to the finite type
stable semistar o\-pe\-rations ${\varphi(\eta)}$ and ${\eta_P}$, defined on $D_P[X]$.
More precisely,
\[
\begin{array}{ll}
 \Spec^{\eta_P}(D_P[X])
 &\hskip -6pt=\{Q\in\Spec(D_P[X])\,|\;Q\cap{N(\star_P)}=\emptyset\}=\\
 &\hskip -6pt=\{Q\in\Spec(D_P[X])\,|\;\mbox{ for all }{g}\in Q\,,\
\content(g)^\star\neq (D_{P})^{\star}\}
\end{array}
\]
and
\[
\begin{array}{ll}
\Spec^{\varphi(\eta)}(D_P[X])
 &\hskip -6pt
=\{Q\in\Spec(D_P[X])\,|\;\varphi^{-1}(Q)\cap{N(\star)}=\emptyset\}=\\
 &\hskip -6pt=\{Q\in\Spec(D_P[X])\,|\;\mbox{ for all
}{f}\in\varphi^{-1}(Q),\content(f)^\star\neq{D^\star}\}=\\
 &\hskip -6pt=\{Q\in\Spec(D_P[X])\,|\;\mbox{ for all
}{f}\in\varphi^{-1}(Q)                                 ,\mbox{ there exists
}\\
 &\hskip -6pt \hspace*{25ex}{H}\in\Spec^\star(D)\mbox{ with
}\content(f)\subseteq{H}\}.
\end{array}
\]
Recall also that the prime ideals in $\Na(D,\star)_{D\setminus{P}} =
D_P[X]^{\varphi(\eta)}$ [respectively, in $\Na(D_P,\star_P) ={D_P[X]^{\eta_P}}$] are in
a na\-tu\-ral bijective correspondence with prime ideals in
$\Spec^{\varphi(\eta)}(D_P[X])$ [respectively, $\Spec^{\eta_P}(D_P[X])$].

\noindent Finally, observe that, in general, $\Spec^{\eta_P}(D_P[X]) \subseteq
\Spec^{\varphi(\eta)}(D_P[X])$, since if $\star_{1} \leq \star_{2}$ are two semistar
operations on an integral domain $R$, then $\Spec^{\star_{2}}(R) \subseteq
\Spec^{\star_{1}}(R)$.

\medskip
\begin{remark} \label{re:3.7} Let $\star$ be a semistar operation
defined on an arbitrary integral domain $D$. Note that: \sl if $ P \in
\Spec^{\star_{f}}(D)$, then $\Na(D,\star)_{D\setminus{P}}=\Na(D_P,\star_P)$. \rm As a
matter of fact, without loss of generality, we can assume that $\star =
\widetilde{\star}$ and, in this case, $PD_{P}$ belongs to
$\Spec^{(\star_{P})_{f}}(D_{P})$ and so \cite[Lemma 2.5 (f)]{fjs}:
$$
\Na(D,\star)_{D\setminus{P}} = \left(\cap \{D_{Q}(X)\mid\, Q\in
\Spec^{\star_{f}}(D)\}\right)_{D\setminus{P}} = D_{P}(X) = \Na(D_P,\star_P)\,.
$$
\end{remark}

\medskip
\begin{proposition} \label{pr:3.8}
 \sl
If $D$ is a B\'ezout domain then, for each $P \in \Spec(D)$ and for each semistar
operation $\star$ on $D$, we have $\Na(D,\star)_{D\setminus{P}}=\Na(D_P,\star_P)$ (and
$\Kr(D,\star)_{D\setminus{P}}=\Kr(D_P,\star_P)$).
\end{proposition}
\begin{proof} If $D$ is a B\'ezout domain, then $\content(g)D$ is a
nonzero principal ideal, for any nonzero polynomial $g\in{D[X]}$. Let
$Q\in\Spec^{\varphi(\eta)}(D_P[X])$, if $Q\nsubseteq{PD_P[X]}$ there exists
$f\in{Q\setminus{PD_P[X]}}$, and, without loss of generality, we may also assume that
$f\in\Ima(\varphi)$ (i.e. $f = \frac{f}{1}$, with $f\in D[X]$), such that
$\content_{D_P}(f)=\content_{D}(f)D_P=D_P$. Henceforth, $\content_{D}(f)=sD$ for some
$s\in{D\setminus{P}}$. Therefore, $\frac{f}{s}\in D[X]$ and
$\content_{D}(\frac{f}{s})=D$. On the other hand, $f\in Q$ and thus
$\frac{f}{s}\in\varphi^{-1}(Q)$, which is a contradiction. As a consequence, we have
$Q\subseteq{PD_P[X]}$ and we conclude that
$\Spec^{\varphi(\eta)}(D_P[X])=\Spec^{\eta_P}(D_P[X])$. As we have noticed above, this
fact implies that $\Na(D,\star)_{D\setminus{P}}=\Na(D_P,\star_P)$. The parenthetical
equality follows from the fact that, for any Pr\"ufer domain $D$ and for any semistar
operation $\star$ on $D$, $\Na(D, \star) = \Kr(D, \star)$ \cite[Remark 3.2]{fjs}, and
so, in particular, $\Na(D_P,\star_P)=\Kr(D_P,\star_P)$, for each $P\in \Spec(D)$.

\end{proof}

Next two examples show that the identity
$\Na(D,\star)_{D\setminus{P}}=\Na(D_P,\star_P)$ does not hold in general.

\begin{example}\label{ex:0208}
 \sl
Let $D:=\mathbb{Z}[Y]$ and let $\star:= d$ be the identical semistar operation on $D$,
i.e., $E^\star:=E$, for any $E\in\oF(D)$. We consider the prime ideal
$P:=2\mathbb{Z}[Y]=2D= (2)$. With the notation introduced in the previous Problem
\ref{pb:3.6}, in case of the local domain $D_{(2)}$, we have:
\[
\begin{array}{ll}
 \Spec^{\eta_{(2)}}(D_{(2)}[X])
 &\hskip
-6pt=\{Q\in\Spec(\mathbb{Z}[Y]_{(2)}[X])\;|\;\,Q\subseteq2\mathbb{Z}[Y]_{(2)}[X]\}\subsetneqq\\

 &\hskip -6pt\subsetneqq\Spec^{\varphi(\eta)}((D_{(2)})[X])=\\
 &\hskip -6pt=\{Q\in\Spec(\mathbb{Z}[Y]_{(2)}[X])\;|\;\mbox{ for all
}{f}\in\varphi^{-1}(Q)\,,\mbox{ there exists}\\
 &\hskip -6pt\hspace*{27ex}{H}\in\Spec(D)\mbox{ with
}\content(f)\subseteq{H}\}.
\end{array}
\]
\rm%
As a matter of fact, if we take $f:=YX+3 \in D[X]$, then
$\content(f)=(Y,3)\neq\mathbb{Z}[Y]=D$, hence there exists a maximal ideal $H$ of $D$
such that $f\in{H[X]}$. In addition, $\varphi(f)$ is not invertible in $D_{(2)}[X]$ and
$\varphi(f)\notin 2D_{(2)}[X]$, hence there exists
$Q\in\Spec^{\varphi(\eta)}(D_{(2)}[X])$, with $Q \neq 2D_{(2)}[X]$, such that
$\varphi(f)\in{Q}$; in particular, $Q\notin\Spec^{\eta_{(2)}}(D_{(2)}[X])$. As a
consequence, we have $\Na(D,\star)_{D\setminus{P}} =\Na(D,
d)_{D\setminus{(2)}}\neq\Na(D_{(2)}, d_{(2)})=\Na(D_P,\star_P)$.

It is also possible to give a direct arithmetic proof of the previous
 fact, that is $\Na(D, d)_{D\setminus{(2)}}\neq \Na(D_{(2)}, d_{(2)})$.  We
consider as before $f:=YX+3\in D[X]$, then $\varphi(f)\in D_{(2)}[X]$ satisfies
$\content_{D_{(2)}}(\varphi(f))=(Y,3)D_{(2)}= D_{(2)}$, hence $\varphi(f)$ is
invertible in $\Na(D_{(2)}, d_{(2)})$.  If $f$ is invertible in $\Na(D,
d)_{D\setminus{(2)}}$, then there exist $h\in D[X]$, $k\in{N(d)}$ and $b\in
D\setminus{(2)}$ such that $fh=kb$.  Since $D$ is a factorial domain, let
$b=p_1p_2\cdots{p_t}$ be a factorization in prime elements of $b$ in $D$.  For any
$p_i$ we have either $p_i\mid{f}$ or $p_i\mid{h}$, hence we may find an identity of the
following type: $fh'=kp_1p_2\cdots{p_s}$, where $h'\in D[X]$, $s\leq t$,
$p_i\nmid{h'}$, for any $i=1, 2, \ldots,s$.

Case 1. If $\content(h')=D$, then
$\content(f)=\content(fh')=\content(k)p_1p_2\cdots{p_s}= p_1p_2\cdots{p_s}$, which is a
contradiction, since $\content(f)$ is not a principal ideal of $D$.

Case 2. If $\content(h')\neq{D}$, then there exists a maximal ideal $H$ in $D$ such
that $\content(h')\subseteq{H}$. From $fh'=kp_1p_2\cdots{p_s}$, with $p_i\nmid{h'}$,
for any $i=1, 2, \ldots,s$, we have $p_1\cdots{p_s}f'h'=kp_1\cdots{p_s}$, for some
$f'\in D[X]$. Therefore, $\content(f'h')= \content(k)=D$. On the other hand,
$\content(f'h') \subseteq\content(h')\subseteq{H}\neq{D}$, which is a contradiction.

From the previous argument we deduce that $f$ is not invertible in  the ring $\Na(D,
d)_{D\setminus{(2)}}$, and so $\Na(D, d)_{D\setminus{(2)}}\neq\Na(D_{(2)}, d_{(2)})$.
\end{example}

\begin{example}
\sl%
Let $D$ be an integral domain and assume that
$D$ possesses two incomparable prime ideals $P_1$ and $P_2$. Set $P:=P_1$ and $\star:=
\star_{\{D_{P_2}\}}$ (i.e. $E^\star := ED_{P_{2}}$, for each $E\in
\boldsymbol{\overline{F}}(D)$), then we have:
\[
\begin{array}{ll}
 \Na(D,\star)_{D\setminus{P}}
 =\left(D[X]_{N(\star)}\right)_{D\setminus{P_1}}
 =\left(D_{P_1}[X]\right)_{N(\star)}\\
 \Na(D_P,\star_P)=D_{P_1}[X]_{N(\star_{P_1})}.
\end{array}
\]
(\,For simplicity of notation, we have identified $D[X]$, and N($\star$), with its
canonical image in $D_{P_1}[X]$\,). We claim that $ D_{P_1}[X]_{N(\star)} \neq
D_{P_1}[X]_{N(\star_{P_1})}$.
\rm%

We compare the two multiplicative sets $N(\star)$ and $N(\star_{P_1})$; more precisely,
we show that there is a canonical injective map
$\left(D_{P_1}[X]\right)_{N(\star)}\longrightarrow{D_{P_1}}[X]_{N(\star_{P_1})}$, which
is not surjective.

\noindent Note that if $f\in N(\star)\ (\subseteq D[X])$, then $\content(f)^\star =
D^\star$, hence $(\content(f)D_{P_1})^{\star_{P_1}} = (\content(f)D_{P_1})^\star
=(D^\star D_{P_1})^\star = D_{P_1}^\star= D_{P_1}^{\star_{P_1}}$, i.e. $f\in
N(\star_{P_1})$.

\noindent Let $f\in{D_{P_1}}[X]$, then there is $s'\in {D\setminus{P_1}}$ such that
$f=s'f'$ with $f'\in{D[X]}$, then --- without loss of generality --- we may assume
$f\in{D[X]}$.

\noindent If $f\in{N(\star_{P_1})}\cap D[X]$, then $\content(f)^{\star_{P_1}}=
(\content(f)D_{P_1})^\star=(D_{P_1})^\star$. On the other hand, the finite type stable
semistar operation $\star$ and the localization at $P_1$ commute, as they are defined
by hereditary torsion theories of finite type. Therefore, we have
$\content(f)^\star{D_{P_1}}=(\content(f)D_{P_1})^\star= (D_{P_1})^\star=
(D^\star)_{D\setminus P_1}$.

\noindent It is obvious that, in general, the previous equality does not imply
$\content(f)^\star=D^\star$, i.e. $f\in N(\star)$. In fact, if
$f\in{P_2[X]}\setminus{P_1[X]}$, then $\content(f)\subseteq{P_2}\setminus{P_1}$ and
$\content(f)D_{P_1}=D_{P_1}$ hence, in particular,
$\content(f)^\star{D_{P_1}}=(D^\star)_{D\setminus P_1}$. On the other hand,
$\content(f)^\star=\content(f)D_{P_2}\subseteq{P_2D_{P_2}}\subsetneqq{D_{P_2}}=D^\star$,
hence $f\notin{N(\star)}$.

\noindent Since $f\in{N(\star_{P_1})}$, then $1/f \in
\left(D_{P_1}[X]\right)_{N(\star_{P_1})}$. If $1/f =h/k$, for some
$h/k\in{D_{P_1}}[X]_{N(\star)}$, where $h \in D_{P_1}[X]$ and $ k\in N(\star)$, then
$hf=k\in{N(\star)}$. Therefore, $\content(hf)^\star=\content(k)^\star=D^\star$, which
is a contradiction, as
$\content(hf)^\star\subseteq\content(f)^\star\subseteq{P_2D_{P_2}}\neq{D_{P_2}}=D^\star$.

\end{example}

Let us now consider the second question considered in Problem \ref{pb:3.6}, concerning
the Kronecker function rings. Also in this case the answer is negative in general, as
the following Example \ref{ex:3.14} proves. First we give a positive example.

\begin{example}
 Let $V$ a valuation overring of an integral domain $D$, with maximal
 ideal $M$, which is not essential (i.e. $D_{P} \subsetneqq V$, where
$P:=M\cap{V}$).
 Set $\star:=\star_{\{V\}}$ (i.e. $E^\star :=
EV$, for each $E\in \boldsymbol{\overline{F}}(D)$). In this situation, $\star$ is a
finite type a.b. semistar operation on $D$ and $P$ is the only maximal
quasi--$\star$--ideal of $D$.

In this particular case, the description of the Kronecker function ring $\Kr(D, \star)$
is rather easy.
Let $a$, $b\in{V}$, we set as usual $\ a\mid_Vb$ \ if there exists $v\in{V}$ such that
$av=b$ and, for each $f\in{D[X]}$, we denote by $a(f)$ a generator (in $V$) of the
principal ideal $\content(f)V$. Then, using also \cite[Example 3.6]{FL1}, we have:
\[
\begin{array}{rl}
 \Kr(D,\star)
 \hskip -5pt &=\left\{\frac{f}{g}\;|\;\,f,g\in{D[X]\setminus\{0\}}\,,\ \mbox{

 with }a(g)\mid_Va(f)\right \} \cup \{0\}=\\
 &= V(X)\,.
\end{array}
\]
Therefore $\Kr(D,\star)_{D\setminus{P}}=\Kr(D,\star)$, because each element $b\in
D\setminus{P}$ is a unit in $V$ (i.e. $\content(b)^\star = (bD)^\star =bV=V$).

 \noindent On the other hand, $\Kr(D_P,\star_P)$ has a similar
 description:
\[
\begin{array}{rl}
 \Kr(D_P,\star_P) \hskip -6pt &=\left\{\frac{bf}{cg}\mid f,
g\in{D[X]\setminus\{0\}},\ b,
 c\in{D\setminus{P}} \mbox{ with } a(cg)\mid_{V} a(bf)\right\}=\\
 \hskip -6pt &= V(X)\,,
 \end{array}
\]
and we conclude immediately that
$\Kr(D_P,\star_P)=\Kr(D,\star)=\Kr(D,\star)_{D\setminus{P}}$.
\end{example}

 Next we give a negative example for the Kronecker function rings case.

\begin{example} \label {ex:3.14}
Let $D:=\mathbb{Z}[\sqrt{-5}]$. Since $D$ is a Dedekind domain, if $d$ is the identical
semistar operation on $D$, then $\Kr(D,d)=\Na(D,d) = D(X)$. We take
$f:=(1+\sqrt{-5})+(1-\sqrt{-5})X \in D[X]$ and we consider, for instance, the prime
ideal $P: =(3, 1+\sqrt{-5})D\nsupseteqq Q:=(1+\sqrt{-5},1-\sqrt{-5})D\ (=\rad(2D))$.
Then, arguing for $f$ as in Example \ref{ex:0208}, we obtain that
$D(X)_{D\setminus{P}}= \Kr(D, d)_{D\setminus{P}}\neq \Kr(D_P,d_P)=D_{P}(X)$.

Note that this example produces also a negative answer to Problem \ref{pb:3.6} (3).
\end{example}

We finish this section with a local characterization of P$\star$MDs. Recall that
P$\star$MDs were characterized in \cite{fjs}, using quasi--$\star_{f}$--prime ideals;
here we extend this characterization by using the whole prime spectrum.  In particular,
next result provides new examples of nontrivial local P$\star$MDs, by taking the
localizations of a P$\star$MD at its prime non quasi--$\star_f$--ideals.

\begin{theorem} \label{cor:3.7}
Let $\star$ be a semistar operation on an integral domain $D$. Then the following
statements are equivalent:
\brom
 \item\sl
 $D$ is a P$\star$MD;
 \item\sl
 $D_{P}$ is a P$\sP$MD, for each $P\in \Spec(D)$;
 \item\sl
 $D_{M}$ is a P$\sM$MD, for each $M\in \Max(D)$;
 \item\sl
 $D_{N}$ is a P$\star_{N}$MD, for each $N\in \calM(\star_{f})$;
 \item\sl $D_{Q}$ is a P$\star_{Q}$MD, for each $Q\in \calQ(\star_{f})$.
 \erom
\end{theorem}

\begin{proof} We already proved that (i) $\Rightarrow$ (ii)
(Remark~\ref{rk:8}). Obviously, (ii) $\Rightarrow$ (iii), (v); and (v) $\Rightarrow$
(iv).

 (iii) (or (ii)) $\Rightarrow$ (i). Recall from \cite[Section 3,
 Theorem 1 and Remark 2]{fjs} that $D$ is a P$\star$MD if and only if
 $\Na(D, \star) = \Kr(D, \star)$ and this last equality follows from
 Proposition~\ref{pr:3.5} (2) and (4).

(iv) $\Rightarrow$ (i). We have already observed (Remark~\ref{re:3.7}) that $\Na(D,
\star)_{D\setminus N} =$ $ \Na(D_{N}, \star_{N})$ and we know that $\Kr(D,
\star)_{D\setminus N} \subseteq \Kr(D_{N}, \star_{N})$, for each $N\in
\calM(\star_{f})$ (Proposition \ref{pr:3.5} (3)). On the other hand, by assumption, $
\Na(D_{N}, \star_{N}) $ $= \Kr(D_{N}, \star_{N})$, for each $N\in \calM(\star_{f})$.
From the previous relations and from the fact that $\Na(D, \star) \subseteq \Kr(D,
\star)$, we deduce immediately that $\Na(D, \star)_{D\setminus N} =\Kr(D,
\star)_{D\setminus N}$, for each $N\in \calM(\star_{f})$ The conclusion follows
immediately, since $ \Na(D, \star) = \cap \{ \Na(D, \star)_{D\setminus N}\,|\; N\in
\calM(\star_{f})\}$ and $ \Kr(D, \star) = \cap \{ \Kr(D, \star)_{D\setminus N}\,|\;
N\in \calM(\star_{f})\}$ (Proposition~\ref{pr:3.5} (2) and (4)).
\end{proof}

\begin{remark}
M. Zafrullah, in \cite{jpaa50}, proves a different local characterization of
P$\star$MDs in the particular case where $\star=v$. More precisely, he obtains that a
domain $D$ is a P$v$MD if and only if (a) $D_P$ is a P$v_{D_P}$MD for every prime ideal
$P$ of $D$ and (b) for every prime $t_D$--ideal $Q$ of $D$, $QD_Q$ is a
$t_{D_Q}$--ideal\ (about condition (b) see also Remark~\ref{newre}). As we have already
observed at the beginning of this section, recall that $v_{D_P}$ is different, in
general, from $v_{P}$.
\end{remark}

\section{Inducing semistar operations}

In this section, we deal with the converse of the problem considered in the first part
of this paper, i.e., we start from a family of ``local'' semistar operations on the
localized rings $D_P$, where $P$ varies in a nonempty set of prime ideals of an
integral domain $D$, and the goal is the description of a gluing process for building a
new ``global'' semistar operation on the ring $D$.

Let $D$ be an integral domain. Let $P$ be a prime ideal of $D$ and let $\ast_P$ be a
semistar operation on the localization $D_P$ of $D$ at $P$. Then we may consider
$\mbox{\d{${\ast}$}$_P$}$, the induced semistar operation on $D$ defined as follows,
for each $E\in\boldsymbol{\overline{F}}(D)$:
\[
 E^{^{\mbox{\tiny \d{${\ast}$}$_P$}}}=(ED_P)^{\ast_P}\;.
\]

Let $\Theta$ be a given nonempty subset of $\Spec(D)$ and let $\{\ast_P \;|\;\, P\in
\Theta \}$ be a family of semistar operations, where $\ast_P$ is a semistar operation
on the localization $D_{P}$ of $D$ at $P$.  We define
 $\boldsymbol{\wedge}:=
 \wedge_{{\Theta},\{\ast_P \}}:=
 \wedge_{{\Theta}}:= \wedge\{{\mbox{\d{${\ast}$}}}_{P} \;| \;\, P\in \Theta \}$
as the semistar operation on $D$ defined as follows, for each
$E\in\boldsymbol{\overline{F}}(D)$,
$$
E^{\boldsymbol{\wedge}}:= \cap \{ (ED_{P})^{\ast_{P}}\;|\;\, P \in \Theta \}.
$$
If $\Theta$ is the empty set, then we set $\boldsymbol{\wedge}:=\wedge_{\emptyset}:=
e_{D}$. Given a semistar operation $\star$ on $D$, for each prime ideal $P$ of $D$, we
denote as usual by $\sP$ the semistar operation $\sdDP$ on $D_{P}$, deduced from
$\star$ by ascent to $D_{P}$ (i.e. $E^{\star_{P}}: = E^\star$, for each $E \in
\overline{\boldsymbol{F}}(D_{P})\ (\supseteq \overline{\boldsymbol{F}}(D)) $; in
particular if $\star $ coincides with the semistar operation
$\boldsymbol{\wedge}:=\wedge_{{{\Theta},\{\ast_P \}}}$ defined on $D$, we can consider
a semistar operation ${\boldsymbol{\wedge}}_{P}$ on $D_{P}$, for each $P\in \Theta$.

Note that:
\balf
\bf\item\rm for each $P \in \Theta$, $ (E^{\boldsymbol{\wedge}} D_{P})^{\ast_{P}}
=(ED_{P})^{\ast_{P}}\,; $

\bf\item\rm for each $P \in \Theta$, ${ \boldsymbol{\wedge}}_{P} \leq \ast_{P}$\,;

\bf\item\rm ${\boldsymbol{\wedge}} = \wedge\{ {\mbox{\d{${\wedge}$}}}_{P} \;| \;\, P\in
\Theta \} = \wedge\{{\mbox{\d{${\ast}$}}}_{P} \;| \;\, P\in \Theta \}\,;$

\bf\item\rm for each semistar operation $\star$ on $D$,
$\star\leq{\underset{\star}{\boldsymbol{\wedge}}}
:=\wedge\{{\mbox{\d{${\star}$}}}_{P}\;| \;\, P\in\Spec(D)\}\,.$
\ealf

(a) We remark that: $E^{
\boldsymbol{\wedge}}=\cap\{(ED_{P'})^{\ast_{P'}}\;|\;P'\in\Theta\}
\subseteq(ED_P)^{\ast_P}$, for each $P\in\Theta$.
Therefore, we deduce that:
\[
 (E^{\boldsymbol{\wedge}}{D_P})^{\ast_{P}}
 \subseteq((ED_P)^{\ast_P}D_P)^{\ast_P}
 =(ED_PD_P)^{\ast_P}
 =(ED_P)^{\ast_P}.
\]

The opposite inclusion is trivial.

(b), (c) and (d) are straightforward.

\begin{theorem}
Let $\star$ be a semistar operation on an integral domain $D$. For each $P\in\Spec(D)$,
denote as usual by $\sP$ the semistar operation $\sdDP$ on $D_{P}$, induced from
$\star$ by ascent to $D_{P}$.  Set
$\underset{\star}{\boldsymbol{\wedge}}:=\wedge\{{\mbox{\d{${\star}$}}}_{P}\;|\;\,P\in\Spec(D)\}$.
If $\star$ is a spectral semistar operation on $D$ then
$\star={\underset{\star}{\boldsymbol{\wedge}}}$.
\end{theorem}
\begin{proof}
 Let $\Delta \subseteq \Spec(D)$ be such that $\star =\star_{\Delta}$.
For each $E\in \boldsymbol{\overline{F}}(D)$, then:
$$
\begin{array}{ll}
 E^\star &= (E^\star)^\star= \left(\cap \{ E D_{Q} \;| \;\, Q \in
 \Delta \}\right)^\star \supseteq\cap \{ (ED_{Q})^\star \;|\;\, Q
 \in \Delta \}=\\
 &= \cap \{ (ED_{Q})^{\star_{Q} } \;| \;\, Q \in \Delta \}=
 \cap \{ E^{{\mbox{\tiny \d{${\star}$}$_{Q}$}}}\;| \;\, Q \in \Delta
 \} \supseteq\\
 & \supseteq \cap \{ E^{{\mbox{\tiny \d{${\star}$}$_{P}$}}}\;| \;\, P \in
 \Spec(D) \}= E^{{\underset{\star}{\boldsymbol{\wedge}}}} = \cap
 \{ (ED_{P})^{\star_{P}}\;| \;\, P \in \Spec(D) \}\supseteq E^\star\,.
\end{array}
$$
As for any $P\in\Delta$ and any $E\in\oF(D)$ we have
$(ED_P)^\star:=\cap\{ED_PD_Q\mid\;Q\in\Delta\}\subseteq{ED_PD_P}=ED_P$.
\end{proof}

\begin{lemma} \label{le:4.1}
Let $D$ be an integral domain and let $P$ be a prime ideal of $D$. If $\ast_{P}$ is a
spectral semistar operation on $D_P$, defined by a nonempty subset ${\Delta}_P\subseteq
\Spec(D_P)$, then {$\mbox{\d{${\ast}$}}$}$_P$ is a spectral semistar operation on $D$
defined by the (nonempty) set $\mbox{\d{${{\Delta}}$}$_P$} := \{Q\in
\Spec(D)\;|\;\,\;QD_P \in {\Delta}_P \}$.
\end{lemma}
\begin{proof}
 For each $E\in\boldsymbol{\overline{F}}(D)$:
\[
\begin{array}{rl}
 E^{^{\mbox{\tiny \d{${\ast}$}$_P$}}}
 &=(ED_P)^{\ast_P}
=\cap \{(ED_P)_H\;|\;\, H\in {\Delta}_{P} \} =\\
&=\cap \{ED_Q\;|\;\, QD_{P}=:H \in {{\Delta}}_P \} = \cap \{ED_Q\;|\;\, Q \in
\mbox{\d{${{\Delta}}$}$_P$} \} = E^{^{\star_{\mbox{\tiny\d{$_{{\Delta}}$}}_{P}}}}.
\end{array}
\]
\end{proof}

\begin{corollary} \label{pr:4.2}
Let $D$ be an integral domain and let $\Theta$ be a nonempty subset of $\Spec(D)$. If
$\ast_{P}$ is a spectral semistar operation on $D_P$, defined by a subset
${\Delta}_{P}\subseteq\Spec(D_P)$, for each $P\in\Theta$, and if \
${\boldsymbol{\wedge}}:=\wedge_{{\Theta},\{\ast_P \}}:=
\wedge\{{\mbox{\d{${\ast}$}}}_{P} \;| \;\,$ $ P\in \Theta \}$ is the semistar operation
on $D$ defined as above, then ${\boldsymbol{\wedge}}$ is a spectral semistar operation
on $D$ defined by the subset ${\Delta}:=\cup\{\mbox{\d{${{\Delta}}$}}_P\;|\;\, P \in
\Theta\}\subseteq\Spec(D)$ (i.e. ${\boldsymbol{\wedge}} = \star_{{\Delta}}$).
\end{corollary}
\begin{proof}
This statement is a straightforward consequence of the previous Lemma~\ref{le:4.1}.
\end{proof}

\begin{lemma} \label{le:4.5}
Let $D$ be an integral domain with quotient field $K$. Let $\Theta$ be a given nonempty
subset of $\Spec(D)$ and let $\{\ast_P \;|\;\, P\in \Theta \}$ be a family of semistar
operations, where $\ast_P$ is a semistar operation on the localization $D_{P}$ of $D$
at $P$. Assume that $\ast_P$ is a semistar operation of finite type and that the family
$\{D_{P}\,|\; P\in \Theta \}$ has the finite character (i.e. for each non zero element
$x \in K$, $xD_{P}= D_{P}$ for almost all the $D_{P}$'s). Then the semistar operation
${\boldsymbol{\wedge}} := {{\wedge_{{{\Theta},\{\ast_P \}}}}:=} \wedge\{
{\mbox{\d{${\ast}$}}}_{P} \;| \;\, P\in \Theta \}$ is a finite type semistar operation
on $D$.
\end{lemma}
\begin{proof}
Let $E\in\boldsymbol{\overline{F}}(D)$, recall that $ E^{{\boldsymbol{\wedge}}} := \cap
\{ (ED_{P})^{\ast_{P}}\;|\;\, P \in \Theta \}$. We want to show that if $x \in
E^{{\boldsymbol{\wedge}}}$ then there exists $F\subseteq E$, with $F
\in\boldsymbol{{f}}(D)$, such that $x \in F^{{\boldsymbol{\wedge}}}$. By the finite
character condition, we may assume that $xD_{P}=D_{P}$, for all
$P\in\Theta\setminus\{P_{1}, P_{2}, \ldots, P_{r}\}$. Since $x \in\cap \{
(ED_{P})^{\ast_{P}}\;|\;\, P \in \Theta \}$, by the finiteness condition on the
$\ast_P$'s, we can find $F_{i} \subseteq E$, with $F_{i} \in \boldsymbol{{f}}(D)$, such
that $xD_{P_{i}}\subseteq (F_{i}D_{P_{i}})^{\ast_{P_{i}}}\subseteq
(ED_{P_{i}})^{\ast_{P_{i}}}$, for $1 \leq i \leq r$, and
$(F_{i}D_{P})^{\ast_P}=(D_{P})^{\ast_P}$, for each
$P\in\Theta\setminus\{P_{1},P_{2},\ldots,P_{r}\}$ and each $1\leq{i}\leq{r}$.

\noindent If we set $F:=F_{1}+F_{2}+\cdots+F_{r}$, then we have that $F \subseteq E$,
with $F \in\boldsymbol{{f}}(D)$, such that $FD_{P}=D_{P} =xD_{P}$, for each $P \in
\Theta \setminus \{P_{1}, P_{2}, \ldots, P_{r}\}$, and
$(FD_{P_{i}})^{\ast_{P_{i}}}\supseteq(F_{i}D_{P_{i}})^{\ast_{P_{i}}} \supseteq
xD_{P_{i}}$, for $1 \leq i \leq r$. We conclude that $x \in F^{{\boldsymbol{\wedge}}} =
\cap \{ (FD_{P})^{\ast_{P}}\;|\;\, P \in \Theta \}$.
\end{proof}

\begin{corollary} \label{cor:4.6}
Let $D$ be an integral domain, let $\Theta$ be a nonempty subset of \ $\Spec(D)$ and
let $\{\ast_P \;|\;\, P\in \Theta \}$ be a family of semistar operations, where
$\ast_P$ is a semistar operation on the localization $D_{P}$ of $D$ at $P\in \Theta$.
We can associate to the semistar operation
$\wedge_{\Theta}:=\wedge_{\Theta,\{\ast_P\}}:=
 \wedge\{{\mbox{\d{${\ast}$}}}_{P}\;|\;\,P\in\Theta\}$ (defined on $D$)
two semistar operations (both defined on $D$): $\widetilde{\wedge_{\Theta}}$ and
$\wedge\{\widetilde{{\mbox{\d{${\ast}$}}}_{P}}\;|\;\,P\in\Theta\} =:
{{\wedge_{{{\Theta},\{\widetilde{\ast_P }\}}}}:=} \wedge_{{\widetilde{\Theta}}}$.
Assume that the family $\{D_{P}\,|\; P\in \Theta \}$ has the finite character, then:
$$
\widetilde{\wedge_{\Theta}} = \wedge_{{\widetilde{\Theta}}}.
$$
\end{corollary}
\begin{proof}
By the previous Lemma~\ref{le:4.5}, $\wedge_{{\widetilde{\Theta}}}$ is a finite type
semistar operation on $D$, since $\widetilde{{\mbox{\d{${\ast}$}}}_{P}}$ is a finite
type semistar operation on $D$, because $\widetilde{{\mbox{{${\ast}$}}}_{P}}$ is a
finite type semistar operation on $D_{P}$, for each $P \in \Theta$, and the family
$\{D_{P}\,|\; P\in \Theta \}$ has the finite character.

Note that $\widetilde{{\mbox{{${\ast}$}}}_{P}}$ is a spectral semistar operation on
$D_P$, defined by the subset ${\Delta}_P := \Spec^{(\ast_{P})_{f}}(D_{P})$. In this
situation, we know from Corollary~\ref{pr:4.2} and Lemma~\ref{le:4.1} that
$\wedge_{{\widetilde{\Theta}}}$ is a spectral semistar operation on $D$ defined by the
set ${\Delta} := \cup \{ \mbox{\d{${{\Delta}}$}}_P \;|\;\, P \in \Theta \}\ (\subseteq
\Spec(D))$, i.e. $\wedge_{{\widetilde{\Theta}}} = \star_{{\Delta}}$.  Therefore, we
deduce that $\widetilde{\wedge_{{\widetilde{\Theta}}}} =
\wedge_{{\widetilde{\Theta}}}$, since $\wedge_{{\widetilde{\Theta}}}$ is a finite type
stable semistar operation \cite[Corollary 3.9 (2)]{FH}. On the other hand,
$\wedge_{{\widetilde{\Theta}}} \leq \wedge_{\Theta}$ (because
$\widetilde{{\mbox{\d{${\ast}$}}}_{P}} \leq {\mbox{\d{${\ast}$}}}_{P}$, for each $P\in
\Theta$), thus we have also that
$\wedge_{{\widetilde{\Theta}}}=\widetilde{\wedge_{{\widetilde{\Theta}}}} \leq
\widetilde{\wedge_{\Theta}}$. Moreover, for each $Q \in {\Delta}$, $QD_{P} \subseteq
(QD_{P})^{(\ast_{P})_{f}} \cap D_{P} \neq D_{P}$, for some $P \in \Theta$. Since
$\left({\wedge_{\Theta}}\right)_{f} \leq {\mbox{(\d{${\ast}$}}}_{P})_{f}$, then $Q
\subseteq Q^{(\wedge_{\Theta})_{f}} \cap D \neq D$.  From this fact, we deduce that
$\widetilde{\wedge_{\Theta}} \leq \star_{{\Delta}}=\wedge_{{\widetilde{\Theta}}} $ and
so we conclude that $ \widetilde{\wedge_{\Theta}} = \wedge_{{\widetilde{\Theta}}}.$
\end{proof}

\begin{theorem} \label{pr:4.7}
Let $D$ be an integral domain, let $\Theta$ be a\ nonempty subset of \ $\Spec(D)$ and
let $\{\ast_P \;|\;\, P\in \Theta \}$ be a family of spectral semistar operations,
where $\ast_P$ is a semistar operation on the localization $D_{P}$ of $D$ at $P \in
\Theta$, defined by a subset $\Delta_{P} \subseteq \Spec(D_{P})$. Set
$\mbox{\d{$\Delta$}}_{P} := \{ Q \in \Spec(D) \,|\; QD_{P} \in\Delta_{P} \}$ and set
${\boldsymbol{\wedge}} :={{\wedge_{{{\Theta},\{\ast_P \}}}}:=}\wedge_{\Theta}$. Assume
that the family of spectral semistar operations $\{\ast_P \;|\;\, P\in \Theta \}$
satisfies the following condition:

{\bf ($\boldsymbol{\downarrow}$)} \ %
for each pair of prime ideals $P$, $P'\in\Theta$, with $P'\neq{P}$, then \\
 {\centerline {$\mbox{\d{$\Delta$}}_{P'}\cap
P^{\downarrow}\subseteq\mbox{\d{$\Delta$}}_{P}$\,.}}

\noindent Set $\Delta:=\{Q\in\Spec(D)\;|\;\,QD_{P}\in\Delta_{P},\,\mbox{ for some
}P\in\Theta\}$. Then, the spectral semistar operation $\star:=\star_{\Delta}$ on $D$
verifies the following properties:
\balf
\item \sl
for each $P \in \Theta$, $\star_{P}=\ast_{P}$ \ (where, as
usual, $\star_{P}:=\dot{\star}^{D_{P}}$);
\item \sl
$(\star_{\Delta}=) \ \star={\boldsymbol{\wedge}}\ (=\wedge_{\Theta})$\, (hence, in
particular, ${\boldsymbol{\wedge}_{P}} = \star_{P}= \ast_{P}$, for each $P\in\Theta$).
\ealf
\end{theorem}
\begin{proof}
(a)
 Fix $P \in \Theta$ and, to avoid the trivial case, assume that $P
 \neq (0)$. Set
$$
\begin{array}{rl}
 \Theta_{0}&:=\{ P'\in \Theta \;|\;\, P' \cap P \,
 \mbox{ does not contain a nonzero prime ideal } \}\,,\;\; \; \mbox{\rm and}\\
 \Theta_{1}&:=\{ P''\in \Theta \;|\;\, P'' \cap P \,
 \mbox{ contains a nonzero prime ideal } \}.
 \end{array}
$$
Note that if $P'\in \Theta_{0}$, then $D_PD_{P'}$ coincides necessarily with $K$, the
quotient field of $D$; note also that $P$ belongs to $\Theta_{1}$.

Assume that $P''\in \Theta_{1}$. We know that there is a bijective correspondence
between prime ideal of $D_PD_{P''}$ and the set $\{H\in\Spec(D)\;| \;\,
H\subseteq{P''}\cap{P}\}$, hence $D_PD_{P''}=\cap\{D_H\;| \;\, H\subseteq{P''}\cap{P}
\mbox{ and }H\in\Spec(D)\} \subsetneq K$. Therefore, by assumption, for each $P''\in
\Theta_{1}$, the set $\mbox{\d{$\Delta$}}_{P''}\cap P^{\downarrow} \subseteq
\mbox{\d{$\Delta$}}_{P}$ and so, for each $G \in \boldsymbol{\overline{F}}(D_{P})$,
$(GD_{P''})^{\ast_{P''}} \supseteq (GD_{P})^{\ast_{P}}= G^{\ast_{P}}$. Henceforth, for
each nontrivial $G\in\overline{\boldsymbol F} (D_{P})$, we have
 $K
 \supsetneq\cap \{(GD_{P''})^{\ast_{P''}} \;|\;\, P''\in \Theta_{1}\}
 =(GD_{P})^{\ast_{P}}
 =G^{\ast_{P}}$;
therefore:
\[
\begin{array}{rl}
 G^{\sP} &= G^{\star} = G^{\star_{\Delta}}= \\
 &=\cap \{\, \cap \{ GD_Q \;| \;\, QD_{H} \in \Delta_{H}\} \;\, | \;\,\, H
\in \Theta \}=\\
 &=\cap \{\, \cap \{ GD_Q \;| \;\, Q\in \Delta,\, Q \subseteq H \}
\;\,|\;\,\, H\in \Theta\}= \\
 &=\cap \{ (GD_{H})^{\ast_{H}} \;| \;\, H \in \Theta\} =\\
 &=(\cap\{(GD_{P'})^{\ast_{P'}}\;|\;\,P'\in\Theta_{0}\})\cap(\cap\{(GD_{P''})^{\ast_{P''}}\;|

 \;\,P''\in\Theta_{1}\})=\\
 &=K \cap ( \cap \{ (GD_{P''})^{\ast_{P''}} \;| \;\, P''\in\Theta_{1}\})=\\
 &=G^{\ast_{P}}.
\end{array}
\]

(b) If $E\in \boldsymbol{\overline{F}}(D)$, then
$$
\begin{array}{ll}
 E^{\star_{\Delta}}
 &=\cap \{ \, \cap \{ ED_Q \;| \;\, QD_{P} \in \Delta_{P}\} \; | \;\, P \in
\Theta \} =\\
 &=\cap \{ \, \cap \{ (ED_{P})D_Q \;| \;\, QD_{P} \in \Delta_{P} \} \; | \;\,
P \in \Theta \} =\\
 &=\cap \{ (ED_{P})^{\ast_{P} }\;| \;\, P \in \Theta \} =\\
 &=E^{\boldsymbol{\wedge}}.
\end{array}
$$
\end{proof}

Next example shows that condition (${\downarrow}$) does not hold in general. Later
(Example \ref{ex:4.12}), we will give an example for which condition (${\downarrow}$)
holds.

\begin{example}\sl Let $D$ be an integral domain and let $\Theta$
be a nonempty subset of $\Spec(D)$ and let $\{\ast_P\;|\;\,P\in\Theta\}$ be a family of
semistar operations, where $\ast_P$ is a semistar operation on the localization $D_{P}$
of $D$ at $P\in \Theta$. Let $\Delta_{P}:={\calQ}((\ast_{P})_{f})=
\Spec^{(\ast_{P})_{f}}(D_P)$ be the set of all the quasi--$(\ast_{P})_{f}$--prime
ideals of $D_{P}$, for each $P\in \Theta$.  The family of spectral semistar operations
$\{\widetilde{\ast_P}\;|\;\,P\in\Theta\}$ does not verify condition \rm
(${\downarrow}$).\rm

For instance, let $D$ be a domain with two incomparable prime ideals $P_1$ and $P_2$
containing a common nonzero prime ideal $Q$. Let  $\ast_{P_1}:=d\ (= d_{D_{P_{1}}})$ be
the identical semistar operation on $D_{P_{1}}$, and let $\ast_{P_2}:=e\
(=e_{D_{P_{2}}})$ be the trivial semistar operation on $D_{P_{2}}$.  We have that
$\ast_{P_1}$ and $\ast_{P_2}$ are both finite type stable semistar operations (i.e.
$\ast_{P_1}= \widetilde{\ast_{P_{1}}}$ and $\ast_{P_2}= \widetilde{\ast_{P_{2}}}$ ),
with $\mbox{\d{$\Delta$}}_{P_1}= \{P\in \Spec(D)\mid P \subseteq P_{1}\}$ and
$\mbox{\d{$\Delta$}}_{P_2}= \{(0)\}$.  The ideal $Q$ produces a counterexample to
condition (${\downarrow}$)\,.\ Indeed $Q\in \mbox{\d{$\Delta$}}_{P_1}\cap
{P_2^\downarrow}$ and $Q\notin \mbox{\d{$\Delta$}}_{P_2}$.

Moreover, set $\Theta :=\{ P_{1}, P_{2}, Q \}$ and $\Theta' :=\{ P_{1}, P_{2} \}$. Let
$\ast_{P_1}$ and $\ast_{P_2}$ be as above and let $\ast_{Q}:= e\ (=e_{D_{Q}})$\ (thus
$\ast_{Q}= \widetilde{\ast_{Q}}$ is also a finite type stable semistar operation). Note
that $\Delta_{P_{1}} = \{PD_{P_{1}}\in \Spec(D_{P_{1}})\mid P \subseteq P_{1}\}$,
$\Delta_{P_{2}} = \{(0)\in \Spec(D_{P_{2}})\}$ and $\Delta_{Q} = \{(0)\in
\Spec(D_{Q})\}$.  In this situation, $\Delta := \mbox{\d{$\Delta$}}_{P_1} \cup
\mbox{\d{$\Delta$}}_{P_2} \cup \mbox{\d{$\Delta$}}_{Q} =\{P\in \Spec(D)\mid P \subseteq
P_{1}\}$.  Therefore, it is easy to see that $\wedge_{\Theta}=\wedge_{\Theta'} $ and it
coincides with the finite type spectral semistar operation $\star:=\star_{\Delta}$, but
$(\wedge_{\Theta})_{Q} =\star_{Q}=d_{Q}\lneq e_{Q} =\ast_{Q}$.
\end{example}

\begin{theorem}
Let $D$ be an integral domain, let $\Theta$ be a nonempty subset of\ $\Spec(D)$ and let
$\{\ast_P \;|\;\, P\in \Theta \}$ be a family of spectral semistar operations, where
$\ast_P$ is a semistar operation on the localization $D_{P}$ of $D$ at $P\in \Theta$,
defined by a subset $\Delta_{P} \subseteq \Spec(D_{P})$. Set\
${\boldsymbol{\wedge}}:=\wedge_{{{\Theta},\{\ast_P \}}}$.  Assume that $\{\ast_P
\;|\;\, P\in \Theta \}$ satisfies the condition\ {\bf ($\boldsymbol{\downarrow}$)}\ and
that\ $\ast_P$ is an e.a.b. [respectively, a.b.] semistar operation on $D_{P}$.  Then
the spectral semistar operation ${\boldsymbol{\wedge}}$ (Theorem~\ref{pr:4.7} (b)) is
also an e.a.b. [respectively, a.b.] semistar operation on $D$.
\end{theorem}
\begin{proof}
Note that from the previous Theorem \ref{pr:4.7}, we have that
${{\boldsymbol{\wedge}}}_{P }= {\ast_{P}}$, for each $P \in \Theta$. Let $F$, $G$,
$H\in\boldsymbol{{f}}(D)$ and suppose that $(FG)^{{\boldsymbol{\wedge}}} \subseteq
(FH)^{{\boldsymbol{\wedge}}} $.  Then, for each $P \in \Theta$, we have
$(FD_{P}GD_{P})^{\ast_{P}}= (FGD_{P})^{\ast_{P}} =
((FG)^{{\boldsymbol{\wedge}}}D_{P})^{\ast_{P}} \subseteq
((FH)^{{\boldsymbol{\wedge}}}D_{P})^{\ast_{P}} =$ $(FD_{P}HD_{P})^{\ast_{P}}$.
Therefore, for each $P \in \Theta$, from the e.a.b. hypothesis on $\ast_{P}$ we have
$(GD_{P})^{\ast_{P}} \subseteq (HD_{P})^{\ast_{P}} $.  We conclude immediately, since
we have that $ G^{{\boldsymbol{\wedge}}}=\cap\{(GD_{P})^{\ast_{P} }\;| \;\, P \in
\Theta \} \subseteq \cap \{ (HD_{P})^{\ast_{P} }\;| \;\, P \in \Theta \}=
H^{{\boldsymbol{\wedge}}}$.  A similar argument shows the a.b. case.
\end{proof}

We apply next the previous theory to the case of the finite type stable (semi)star
operation $w := \widetilde{v}$ canonically associated to the (semi)star operation $v$.

\begin{corollary} \label{co:4.3}
 Let $D$ be an integral domain. For each $P \in
 \Spec(D)$, let $w_{D_P}:= \widetilde{v_{D_{P}}}$ be the
 finite type spectral (semi)star operation on $D_P$, defined by the set\ $
 \Spec^{t_{D_{P}}}(D_P)$ of all the $t$--prime ideals of $D_{P}$.  If
 $$
 \underset{w}{\boldsymbol{\bar{\wedge}}}
 :=\wedge\{{\mbox{\d{${w}$}}}_{D_P}\;|\;\,P\in\Spec(D)\}
 $$
 then $\underset{w}{\boldsymbol{\bar{\wedge}}}$ is a
 spectral (semi)star operation on $D$ defined by the following set of
 prime ideals of $D$:
 $$
 \Upsilon:=\cup\mbox{\Large{\{}}\,\{Q\in\Spec(D)\;|\;\,
 QD_{P}\in\Spec^{t_{D_P}}(D_P)\;\}\;\,|\;\;P\in\Spec(D)\mbox{\Large{\}}}\,,
 $$
 i.e.
 $\underset{w}{\boldsymbol{\bar{\wedge}}}=\star_{\Upsilon}$.
\end{corollary}
\begin{proof}
This statement is a particular case of Corollary~\ref{pr:4.2}.
\end{proof}

At this point, it is natural to investigate the relationship between the spectral
(semi)star operation $\underset{w}{\boldsymbol{\bar{\wedge}}}$ (considered in the
previous Corollary \ref{co:4.3}) and the finite type spectral (semi)star operation,
 $w_{D}: = \widetilde{v_{D}}$, on $D$ defined by the set
 $\Spec^{t_{D}}(D)$ of all the $t$--prime ideals of $D$.
We will see that, in general, they are different.

\begin{lemma}\label{le:020712}
Let $D$ be an integral domain. For any prime ideal $P$ of $D$, we denote by $t_{P}$ the
semistar operation $(t_{D})_{P}$ of $D_{P}$\ (defined by $ E^{t_{P}} := E^{t_{D}} =
\cup \{F^{v_{D}}=(D : (D : F)) \mid F \subseteq E \mbox{ and } F \in
\boldsymbol{f}(D)\}$, for each $ E \in \overline{\boldsymbol{F}}(D_{P})$). Let $P$ be a
prime ideal of $D$, we have that $P\in\Spec^{t_{D}}(D)$ if and only if
$PD_P\in\Spec^{t_{P}}(D_P)$.
\newline
In addition, for any prime ideal $P$,  we have
$\Spec^{t_{D_P}}(D_P)\subseteq\Spec^{t_P}(D)$.
\end{lemma}
\begin{proof}
($\Rightarrow$).  %
Assume that, for some $P\in\Spec^{t_{D}}(D)$, we have $(PD_{P})^{t_P} \supsetneq
PD_{P}$.  Then there exists $z \in (PD_{P})^{t_P} = (PD_{P})^{t_D} = \cup \{F^{v_{D}}
\mid F \subseteq PD_{P} \mbox{ and } F \in \boldsymbol{f}(D) \}$, but $z \notin
PD_{P}$. Hence, for some $F \subseteq PD_{P} \mbox{ and } F \in \boldsymbol{f}(D)$, $z
\in F^{v_{D}} \setminus PD_{P}$.  Since $F$ is finitely generated and $F \subseteq
PD_{P}$ then, for some $b\in D \setminus P$, we have that $bF \subseteq PD_{P} \cap D =
P$.  Therefore, $bz \in bF^{v_{D}} = (bF)^{v_{D}} \subseteq P^{t_{D}} =P$, thus $z \in
b^{-1}P \subseteq PD_{P}$, which is a contradiction.

($\Leftarrow$). %
Assume that $(PD_{P})^{t_P}= PD_{P}$.  Note that $PD_{P} = (PD_{P})^{t_P} =
(PD_{P})^{t_D}$ $\supseteq P^{t_D} $.  Henceforth, $P = PD_{P} \cap D \supseteq P^{t_D}
\cap D = P^{t_D}$, hence $P = P^{t_D}$.

For the final statement we proceed as follows.  Let $PD_P\in\Spec^{t_{D_P}}(D_P)$ and
let $F$ be a finitely generated ideal of $D$ contained in $P$, then $F^{v_{D}}
\subseteq (FD_P)^{v_{D} }\subseteq (FD_P)^{v_{D_{P}}} \subseteq (PD_P)^{t_{D_{P}}} =
PD_P$.  Therefore $P^{t_{D}}\subseteq PD_P$ and so $P^{t_{D}}=P$.
\end{proof}

\begin{remark}\label{re:020712}
The same proof given above (Lemma \ref{le:020712}) shows the following general
statement:
\sl%
Let $P$, $Q$ be two prime ideals of an integral domain $D$, then
$PD_Q\in\Spec^{t_Q}(D_Q)$, for each prime ideal $Q$, with $P \subseteq Q$, if and only
if $P\in\Spec^{t_{D}}(D)$.
\end{remark}

\begin{remark}\label{newre}
We emphasize that, in general, the semistar operation $t_P$ does not coincide with the
(semi)star operation $t_{D_P}$, i.e. \emph{$t_P$ is not the $t$--operation on $D_P$}.
For a prime $t_{D}$--ideal $P$ of $D$, the question of when the extended ideal $PD_P$
is a $t_{D_P}$--ideal was studied by M. Zafrullah in \cite{jpaa50} and \cite{jpaa65}\
(where the $t_{D}$--primes $P$ of $D$ such that $PD_P$ is a $t_{D_P}$--ideal were
called \emph{ well behaved prime $t$--ideals}).

For instance, if $P$ is not a well behaved prime $t$--ideal of $D$, then necessarily
$PD_{P}= (PD_{P})^{t_P} \subsetneq (PD_{P})^{t_{D_P}}$.

Using the same argument of the proof of the last statement of Lemma \ref{le:020712},
note that, if $PD_Q$ is a $t_{D_Q}$--ideal, for some prime ideal $Q$ containing $P$)
then $P$ is a $t$--ideal of $D$.  Therefore, using  Remark \ref{re:020712}, we have: if
$Q\in\Spec(D)$ satisfies $Q\supseteq{P}$ and $PD_Q\in\Spec^{t_{D_Q}}(D_Q)$, then
$P\in\Spec^{t_{D}}(D)$, and this happens if and only if for any $Q\in\Spec(D)$, such
that $Q\supseteq{P}$, we have $PD_Q\in\Spec^{t_{Q}}(D_Q)$.
\end{remark}

\begin{example} \label{ex:4.12} \sl The set of all the $t$--prime ideals of
an integral domain $D$ induces a ``natural'' example for which condition \rm
(${\downarrow}$) \sl of Theorem~\ref{pr:4.7} holds.  \rm

\noindent For each $P\in\Spec(D)$, we consider on $D_P$ the set
$\Omega_{P}:=\Spec^{t_{P}}(D_{P})$.  Let $\omega_{P}$ be the spectral semistar
operation on $D_{P}$, defined by $\Omega_{P}$, i.e. $\omega_{P} :=\star_{\Omega_{P}}$.
From Remark \ref{re:020712}, we deduce immediately that
$\mbox{\d{$\Omega$}}_{P'}\cap{P^{\downarrow}}\subseteq \mbox{\d{$\Omega$}}_{P}$, for
each pair $P$, $P'\in\Spec(D)$ such that $P\neq{P'}$. Therefore, the family of spectral
semistar operations $\{\omega_{P} \mid P\in \Spec(D)\}$ verifies condition
(${\downarrow}$).
\end{example}

\begin{corollary} \label{co:4.13}
 Let $D$ be an integral domain.  Let $w_{D}:= \widetilde{v_{D}}$ be
 the finite type spectral (semi)star operation on $D$, defined by the set\ $
 \Spec^{t_{D}}(D)$ of all the $t$--prime ideals of $D$.  For each $P \in
 \Spec(D)$, set as usual $w_{P}:=(w_{D})_{P}$ and let $ w_{D_{P}} :=
 \widetilde{v_{D_{P}}}$ [respectively, $ \omega_{P}$] be the spectral
 semistar operation on $D_P$, defined by the set \ $
 \Spec^{t_{D_{P}}}(D_P)$ [respectively, $\Spec^{t_{{P}}}(D_P)$].  Then:
 $$
 \begin{array}{rl}
 w_{D} =
 & \hskip -4pt \wedge\{ {\mbox{\d{${\omega}$}}}_{P}\;| \;\, P\in \Spec^{t_{D}}(D) \}
 =\wedge\{{\mbox{\d{${w}$}}}_{P} \;| \;\, P\in\Spec^{t_{D}}(D)\}= \\
 =& \hskip -4pt \underset{w}{\boldsymbol{\wedge}}\  (:=
 \wedge\{{\mbox{\d{${w}$}}}_{P} \;| \;\, P\in \Spec(D) \}) \leq
 \underset{w}{\boldsymbol{\bar{\wedge}}}\
 (:=\wedge\{{\mbox{\d{${w}$}}}_{D_P}\;|\;\,P\in\Spec(D)\})\,.
\end{array}
$$
\end{corollary}
\begin{proof}
From Theorem \ref{pr:4.7} (and Example \ref{ex:4.12}), we have that
$\wedge\{ {\mbox{\d{${\omega}$}}}_{P} \;| \;\, P\in \Spec^{t_{D}}(D) \}$ is the
spectral semistar operation $\star_{\Omega}$, where $\Omega := \{Q \in \Spec(D) \mid
QD_{P}\in \Spec^{t_{{P}}}(D_P), \mbox{ for some } P \in \Spec^{t_{D}}(D) \}$. It is
easy to see that $\Omega = \Spec^{t_{D}}(D)$, hence $w_{D} = \star_{\Omega} = \wedge\{
{\mbox{\d{${\omega}$}}}_{P} \;| \;\, P\in \Spec^{t_{D}}(D) \} $.  Moreover, again from
Theorem \ref{pr:4.7}, we have that $ w_{P} = (\star_{\Omega})_{P} = \omega_{P}$, for
each $P\in \Spec^{t_{D}}(D)$, hence $\wedge\{ {\mbox{\d{${w}$}}}_{{P}} \;| \;\, P\in
\Spec^{t_{D}}(D) \} = \wedge\{ {\mbox{\d{${\omega}$}}}_{P} \;| \;\, P\in
\Spec^{t_{D}}(D) \} = w_{D}$.  The last inequality in the statement is a consequence of
Corollary \ref{co:4.3}, since by Lemma \ref{le:020712} $\Upsilon \subseteq \Omega $ and
thus $w_{D} =\star_{\Omega}\leq \star_{\Upsilon } =
\underset{w}{\boldsymbol{\bar{\wedge}}} $.
\end{proof}

\end{document}